\documentclass[11pt]{article}
\usepackage{amsmath, amssymb, color, amsthm,graphicx}
\usepackage{hyperref, enumerate}

\numberwithin{equation}{section}

\newenvironment{proofof}[1]{\begin{proof}[Proof of #1]}{\end{proof}}

\newcommand{\be}{\begin{equation}}
\newcommand{\ee}{\end{equation}}
\newcommand{\bes}{\begin{equation*}}
\newcommand{\ees}{\end{equation*}}
\newcommand{\ds}{\displaystyle}

\newcommand{\scs}{\scriptstyle}

\newcommand{\Det}[2]{\overset{#2}{\underset{#1}{\det}}}

\newcommand{\weighta}{a}
\newcommand{\weightb}{b}

\title{New enumeration formulas for alternating sign matrices and square ice partition functions}
\author{Arvind Ayyer\footnote{Supported by the National Science Foundation under grant DMS-0955584.} \and Dan Romik\footnote{Supported by the National Science Foundation under grant DMS-0955584.}}

\newcommand{\demph}{\textbf}
\newcommand{\asm}{\operatorname{ASM}}

\newtheorem{thm}{Theorem}
\newtheorem{lem}[thm]{Lemma} 
\newtheorem{cor}[thm]{Corollary} 
\newtheorem{prop}[thm]{Proposition}

\begin{document}

\maketitle

\begin{abstract}
The refined enumeration of alternating sign matrices (ASMs) of given order having prescribed behavior near one or more of their boundary edges has been the subject of extensive study, starting with the Refined Alternating Sign Matrix Conjecture of Mills-Robbins-Rumsey \cite{millsetal2}, its proof by Zeilberger \cite{zeilberger2}, and more recent work on doubly-refined and triply-refined enumeration by several authors. In this paper we extend the previously known results on this problem by deriving explicit enumeration formulas for the ``top-left-bottom'' (triply-refined) and ``top-left-bottom-right'' (quadruply-refined) enumerations. The latter case solves the problem of computing the full boundary correlation function for ASMs. 
The enumeration formulas are proved by deriving 
new representations, which are of independent interest,
for the partition function of the square ice model with domain wall boundary conditions at the ``combinatorial point'' $\eta={2\pi/3}$.
\end{abstract}

\section{Introduction} \label{sec:intro}

\subsection{Alternating sign matrices}

An \demph{alternating sign matrix (ASM)} of order $n$ is an $n\times n$ square matrix with entries in $\{0,1,-1\}$ such that in every row and every column, the sum of the entries is $1$ and the non-zero terms appear with alternating signs; see Figure~\ref{fig-asm-example} for an example. From this seemingly innocuous definition a uniquely fascinating class of objects arises: originally discovered by Mills, Robbins and Rumsey in connection with their study of Dodgson's condensation method for computing determinants, ASMs have since been found to have deep connections to many other topics of interest in combinatorics and statistical physics. Some places where ASMs make an unexpected appearance are: the \demph{square ice model} (a.k.a.\ the \demph{six-vertex model}) \cite{karklinskyromik, kuperberg, stroganov, zeilberger2}; \demph{totally symmetric self complementary plane partitions} \cite{fonsecaetal, millsetal3}; \demph{descending plane partitions} \cite{millsetal2}; \demph{domino tilings}~\cite{aztecdiamond}; and the \demph{$O(1)$ loop model} in a cylindrical geometry \cite{cantinisportiello, razumovstroganov}.

\begin{figure}
\begin{center}
$  \left( \begin{array}{cccccc}
0 & 0 & 1 & 0 & 0 & 0 \\
0 & 1 & -1 & 0 & 1 & 0 \\
1 & -1 & 0 & 1 & 0 & 0 \\
0 & 1 & 0 & 0 & -1 & 1  \\
0 & 0 & 1 & 0 & 0 & 0 \\
0 & 0 & 0 & 0 & 1 & 0 
\end{array} \right)  
$
\caption{An alternating sign matrix of order 6.
\label{fig-asm-example}}
\end{center}
\end{figure}

\subsection{Enumeration of alternating sign matrices}

The current paper will be focused on one particular aspect of the study of ASMs, namely the problem of enumerating various naturally-occurring sets of ASMs of some fixed order $n$. This is a problem with a venerable history, starting with the seminal paper \cite{millsetal2} of Mills, Robbins and Rumsey. Having defined ASMs and discovered the role that they play in the definition of the \demph{$\lambda$-determinant}, a natural generalization of matrix determinants, Mills et\ al.\ considered the problem of finding the total number $A_n$ of ASMs of order~$n$. Based on numerical observations, they conjectured the formula
\begin{equation} \label{eq:asm-thm}
A_n = \frac{1!\,4!\,7!\cdots (3n-2)!}{n!(n+1)!\cdots(2n-1)!} = \prod_{j=0}^{n-1} \frac{(3j+1)!}{(n+j)!}
= \prod_{j=0}^{n-1} \frac{\binom{3j+1}{j}}{\binom{2j}{j}}.
\end{equation}
Another natural enumeration problem concerned the so-called \emph{refined} enumeration of ASMs. It is based on the trivial observation (an immediate consequence of the definition of ASMs) that the top row of an ASM contains a single $1$ and no $-1$s. The position of the $1$ in the top row is therefore an interesting parameter by which one  may refine the total enumeration. Thus, for $1\le k\le n$ Mills\ et\ al. defined
$$ A_{n,k} = \#\textrm{ of ASMs of order $n$ with $1$ in position $(1,k)$},
$$
and conjectured that
\begin{equation} \label{eq:refined-asm-thm}
A_{n,k} = \binom{n+k-2}{k-1} \frac{(2n-k-1)!}{(n-k)!} \prod_{j=0}^{n-2} \frac{(3j+1)!}{(n+j)!}  \qquad (1\le k\le n).
\end{equation}
The conjectures \eqref{eq:asm-thm} and \eqref{eq:refined-asm-thm} became the subject of intensive study by combinatorialists. The former became known as the Alternating Sign Matrix Conjecture, and the latter as the Refined Alternating Sign Matrix (RASM) Conjecture; both were eventually proved by Zeilberger \cite{zeilberger1, zeilberger2} (in the case of \eqref{eq:refined-asm-thm}, building on techniques introduced by Kuperberg in \cite{kuperberg}). A readable account of these developments can be found in the book by Bressoud \cite{bressoud}.

\subsection{Doubly-refined enumeration and beyond} 
By symmetry, the observation mentioned above concerning the behavior of the top row of an ASM applies not just to the top row but also to the bottom row and to the leftmost and rightmost columns. This lead Mills et al.\ to consider a \emph{doubly-refined} enumeration involving two parameters, one for the position of the $1$ in the top row and another for the corresponding position in the bottom row. Thus, for $1\le i,j\le n$ we may denote
$$ A_n^{\textrm{TB}}(i,j) = \#\textrm{ of ASMs of order $n$ with $1$ in positions $(1,i)$, $(n,j)$}.
$$
Mills et al.\ did not conjecture an explicit formula for $A_n^{\textrm{TB}}(i,j)$. However, they discovered that a different
family of objects, namely the so-called totally symmetric self-complementary plane partitions (TSSCPPs), had a two-parameter refinement which they conjectured \cite{millsetal3} is also given by the same family of numbers $\big(A_n^{\textrm{TB}}(i,j)\big)_{1\le i,j\le n}$. This was proved in recent years by Fonseca and Zinn-Justin \cite{fonsecaetal}. Also fairly recently, Stroganov \cite{stroganov} derived an explicit formula for these ``top-bottom'' enumeration coefficients, expressing them in terms of the singly-refined coefficients $A_{n,k}$. His result states that, for $1\le i\le j\le n$, 
\begin{equation} \label{eq:stroganov}
A_n^{\textrm{TB}}(i,j) = A_{n,j-i} + \sum_{k=1}^{i-1} D_n(k,j-i+k), 
\end{equation}
where
$$ D_n(s,t) = \frac{1}{A_{n-1}}\left( A_{n-1,t} (A_{n,s+1}-A_{n,s}) + A_{n-1,s} ( A_{n,t+1}-A_{n,t}) \right). $$
(To be a bit more precise, Stroganov proved that the numbers $A_n^{\textrm{TB}}(i,j)$ satisfy the recurrence relation $A_n^{\textrm{TB}}(i+1,j+1)-A_n^{\textrm{TB}}(i,j)=D_n(i,j)$ for all $1\le i,j\le n$; this immediately implies \eqref{eq:stroganov} by summation.)
Stroganov also considered another natural doubly-refined enumeration, the ``top-left'' enumeration of ASMs with prescribed behavior in the top row and left column. For \hbox{$1\le i,j\le n$} denote as above
$$ A_n^{\textrm{TL}}(i,j) = \#\textrm{ of ASMs of order $n$ with $1$ in positions $(1,i)$, $(j,1)$}.
$$
Stroganov proved that the top-left coefficients $(A_n^{\textrm{TL}}(i,j))_{i,j}$ are related to the top-bottom coefficients $(A_n^{\textrm{TB}}(i,j))_{i,j}$ via the relations
\be
\begin{split}
\label{eq:stroganov2}
A_n^{\textrm{TB}}(i,j) &= 
A_n^{\textrm{TL}}(i,j+1) + A_n^{\textrm{TL}}(i+1,j) - A_n^{\textrm{TL}}(i+1,j+1)\ \ \ \ (i,j\ge 2),
\\
A_n^{\textrm{TL}}(2,2) & = A_{n-1},
\end{split}
\ee
which were also rederived by Fischer \cite{fischer4} using different methods.
It is not difficult to invert these relations and therefore to obtain an explicit formula for $A_n^{\textrm{TL}}(i,j)$ (see \cite[Section 3]{fischer4}):
\begin{align*}
A_n^{\textrm{TL}}(i,j) &=
\begin{cases} 
A_{n-1} & 
\textrm{if }i=j=1, \\ 
0 & 
\textrm{if }
i=1<j \\[-5pt]
    &  \ \ \ \textrm{or }j=1<i, \\[-3pt]
\binom{i+j-4}{i-2}A_{n-1}  - \sum_{p=1}^{i-1}\limits
\sum_{q=1}^{j-1}\limits \binom{i+j-2-p-q}{i-1-p} A^{\textrm{TB}}_{n}(p,q)
&
\textrm{if }
i,j\ge 2.
\end{cases}
\end{align*}

A further development of the theory of refined enumerations of ASMs came with the papers \cite{fischerromik}, \cite{karklinskyromik}. In \cite{fischerromik}, Fischer and Romik introduced a new family of doubly-refined enumeration coefficients, which we will call here the ``top-two'' (or ``top-top'') coefficients, since they enumerate ASMs based on their behavior in the top two rows. 
The definition of these coefficients is based on the slightly subtle observation that for any ASM $M=(m_{i,j})_{i,j=1}^n$ of order $n$, there is a unique pair $(i,j)$ with $1\le i<j\le n$ such that $m_{1,i}+m_{2,i}=m_{1,j}+m_{2,j}=1$, and that the top rows of $M$ are uniquely determined by specifying the pair $(i,j)$ along with the unique $k$ such that $m_{k,1}=1$, which satisfies $i\le k\le j$ (thus for a given pair $(i,j)$ there will be $j-i+1$ possible values for $k$). So it makes sense to define the top-two enumeration coefficients as the integers
\begin{align*}
&A_n^{\textrm{TT}}(i,j) = \frac{1}{j-i+1} \\ &
\ \ \ \ \ \times\,\#\Big\{ M=(m_{i,j})_{i,j=1}^n\in\asm_n \ |\ m_{1,i}+m_{2,i} = m_{1,j}+m_{2,j} = 1 \Big\},
\end{align*}
where $1\le i<j\le n$ and $\asm_n$ denotes the set of ASMs of order $n$. (The more correct way of thinking about this enumeration is in terms of complete monotone triangles, a class of objects that is in bijection with ASMs; see \cite{fischerromik} for more details.)
Fischer and Romik studied the problem of finding a formula for $A_n^{\textrm{TT}}(i,j)$. They managed to derive a system of linear equations satisfied for each $n$ by the coefficients $(A_n^{\textrm{TT}}(i,j))_{i,j}$, which, modulo a reasonable conjecture on the invertibility of the system, could be used to express them explicitly as ratios of determinants by Cramer's rule. They also conjectured a more explicit but very complicated formula expressing $A_n^{\textrm{TT}}(i,j)$ in terms of a hypergeometric summation. Finally, Karklinsky and Romik \cite{karklinskyromik} and, shortly afterwards, Fischer \cite{fischer3}, derived two much simpler explicit formulas. According to Karklinsky and Romik's formula, for $1\le i<j\le n$, $A_n^{\textrm{TT}}(i,j)$ is given by
\be \label{eq:toptop-karkromik}
A_n^{\textrm{TT}}(i,j) = \sum_{p=0}^{n-j} \sum_{q=0}^p (-1)^q \binom{p}{q} E_n(i+q,j+p),
\ee
where
$$ 
E_n(s,t) = \frac{1}{A_{n-1}}\left( A_{n-1,t} (A_{n,s+1}-A_{n,s}) - A_{n-1,s} ( A_{n,t+1}-A_{n,t}) \right).
$$
(Note the similarity in the definitions of $D_n(s,t)$ and $E_n(s,t)$, which differ only by a single sign; it is also easy to check that $E_n(s,t) = D_n(s,n-t)$ by using the symmetry $A_{n,k} = A_{n,n+1-k}$.) Fischer's formula on the other hand states that, in our current notation,
\be
\label{eq:toptop-fischer}
A_n^{\textrm{TT}}(i,j) = \sum_{k=j}^n (-1)^{n+k} \binom{2n-2-j}{k-j}
A_n^{\textrm{TB}}(i,k).
\ee
We are not aware of a direct way to establish the
equivalence of the two formulas \eqref{eq:toptop-karkromik} and \eqref{eq:toptop-fischer}.

Going beyond the doubly refined enumeration, Fischer and Romik started developing a theory for the ``$k$-tuply refined'' enumeration of ASMs with prescribed behavior in the top $k$ rows. For each $k\ge 1$, they defined a family of coefficients $(A_n^{\textrm{T}^k}(j_1,\ldots,j_k))_{1\le j_1<\ldots< j_k\le n}$ that encode this enumeration (the correct way to do it is to think about ASMs as \demph{complete monotone triangles}; 
see \cite{fischerromik}),
 and conjectured a generalization of the system of linear equations proved to hold for the top-two doubly-refined coefficients. That conjecture was later proved by Fischer \cite{fischer3}. Recently, Fischer also derived additional linear equations relating the ``top-left-bottom'' and ``top-top-bottom'' families of triply-refined enumeration coefficients \cite{fischer4}. However, these results left open the problem of deriving explicit, closed-form formulas to compute these coefficients.

\subsection{New results: triply- and quadruply-refined enumeration}

Our goal in this paper is to extend the known results on enumeration of alternating sign matrices with prescribed behavior at one or more of their boundary edges, to fully take into account the joint behavior with respect to all edges, thus effectively completing the study of this type of boundary enumeration.

We 
consider the triply-refined and quadruply-refined enumerations for ASMs with prescribed behavior in three (resp.\ four) of their boundary rows/columns. The triply refined enumeration coefficients, which we will refer to as ``top-left-bottom'' coefficients, are given for any $1\le i,j,k \le n$ by
$$
A_n^{\textrm{TLB}}(i,j,k) = \#\Big\{ M=(m_{i,j})_{i,j=1}^n \in \asm_n\ |\ m_{1,i} = m_{j,1} = m_{n,k} = 1 \Big\}.
$$
Similarly, the quadruply-refined ``top-left-bottom-right'' enumeration coefficients are defined for $1\le i,j,k,\ell\le n$ by
$$
A_n^{\textrm{TLBR}}(i,j,k,\ell) = \#\Big\{ (m_{i,j})_{i,j=1}^n \in \asm_n\ |\ m_{1,i}\!=\!m_{j,1}\!=\!m_{n,k}\!=\!m_{\ell,n}\!=\!1 \Big\}.
$$
We remark that in the statistical physics literature, the vector of enumeration coefficients measuring prescribed behavior along specific boundary edges (or in some situations the generating function of this vector) is often referred to as the \demph{boundary correlation function}.
Boundary correlation functions have been studied in great generality by others in the context of integrable systems \cite{bogkitzvo,bogprozvo,fodapres,motegi}.

Our first main result is the following explicit formula relating the top-left-bottom coefficients to the usual singly-refined coefficients $(A_{n,k})_{k=1}^n$. It turns out to be most natural to express the formula as an identity between two multivariate generating functions.

We define the following (single-variable) generating functions $\alpha_n(t), \beta_n(t)$, 
$\gamma_n(t)$ and 
$\delta_{n}(t)$. 
Note that, subsequent to our release of the preprint version of this paper, we were informed by Filippo Colomo \cite{colomo-private} that these formulas can be simplified considerably; see the discussion following Theorem~\ref{thm:quad} below. The generating functions are defined by
\be
\begin{split} \label{defalph}
\alpha_n(t) =& \sum_{k=1}^{n} A_{n,k} t^{k-1},
\\
\beta_n(t) =&  \sum_{k=1}^{n-1} A_{n-1,k} t^{k} = t\, \alpha_{n-1}(t),
\\
\gamma_n(t) =& \sum_{k=1}^{n+2} \Big(
-2(n-k+3) A_{n-1,k-3} + (5n-4k+6) A_{n-1,k-2} \\ & 
\qquad + (n+4k-6) A_{n-1,k-1} -2k A_{n-1,k} 
\Big) \, t^{k-1}, \\
\delta_n(t) =& \sum_{k=1}^{n+3} \Big( 4(n+4-k)(n+5-k) A_{n-1,k-5} \\
&-4(n+4-k)(5n+11-4k) A_{n-1,k-4}  \\
&+(240 - 172 k + 32 k^2 + 120 n - 52 k n + 21 n^2) A_{n-1,k-3} \\
&-2(80 - 80 k +20 k^2 +42 n - 20 k n -5 n^2) A_{n-1,k-2} \\
& +(64 - 84 k + 32 k^2 -4 n -12 k n + n^2) A_{n-1,k-1} \\
&-4 k (n-5+4k) A_{n-1,k}
+4 k (k+1) A_{n-1,k+1} \Big) t^{k-1}.
\end{split}
\ee
Denote by $\Delta(t_1,\ldots,t_k)=\prod_{1\le i<j\le k} (t_i-t_j)$ the standard Vandermonde product of indeterminates $t_1,\ldots,t_k$.

\begin{thm}[Triply refined boundary correlation function] \label{thm:trip}
The generating function of ASMs refined according to the positions of the 1s in the first row, leftmost column and last row satisfies
\be
\begin{split} \label{tripref}
(1-y&+xy)(1-z+yz)
\left( \sum_{i=2}^n \sum_{j=2}^{n-1} \sum_{k=2}^n 
A_n^{\textrm{TLB}}(i,j,k) x^{i-2} y^{n-j-1} z^{n-k} \right)
\\ = & \,
\rho_n \Delta(x,y,z)^{-1}
\det\left( \begin{array}{ccc}
  (x-1)^2 \alpha_n(x) & (y-1)^2\alpha_n(y) & (z-1)^2\alpha_n(z) \\
  (x-1) \beta_n(x) & (y-1)\beta_n(y) & (z-1)\beta_n(z) \\
  \gamma_n(x) & \gamma_n(y) & \gamma_n(z)
\end{array}\right)
\\ &
-(1-y+xy) z^{n-1} \alpha_{n-1}(x)
- (1-z+yz) y^{n-2} \alpha_{n-1}(z),
\end{split}
\ee
where $\rho_n$ is a constant given by
\be \label{defrhon}
\rho_n = \frac{1}{8 A_n^2(2n-2)} \left(\frac{(n-2)!(3n-2)!}{(2n-3)!(2n-1)!}\right)^2.
\ee
\end{thm}

Our next result gives an analogous expansion for the generating
function of the quadruply refined enumeration coefficients
$(A_n^\textrm{TLBR}(i,j,k,\ell))$. Naturally, the expressions become more
complicated, but the result is structurally similar.
It will be convenient to define the doubly refined generating function,
\be \label{deftopleftgf}
\mathcal A_{n}^{\textrm{TL}}(x,y) = \sum_{i,j=2}^{n} A^{TL}_{n}(i,j) x^{i-2} y^{j-2}.
\ee
Note that $\mathcal A_{n}^{\textrm{TL}}(x,y)$ can be expressed up to a constant (that depends on $n$) 
in terms of the functions 
$\alpha_{n}$ and $\beta_{n}$ defined in \eqref{defalph} by combining equations (32) and (33) of \cite{stroganov}.

\begin{thm}[Full boundary correlation function] \label{thm:quad}
The generating function of ASMs refined according to the positions of the 1s in the first row, leftmost column, last row and rightmost column satisfies
\be
\begin{split}  \label{quadref}
&\prod_{t=1}^{4} (1-x_{t+1}+x_{t} x_{t+1})
\left( \sum_{i=2}^{n-1} \sum_{j=2}^{n-1} \sum_{k=2}^{n-1} \sum_{\ell=2}^{n-1} 
A_n^{\textrm{TLBR}}(i,j,k,\ell) x_1^{i-2} x_2^{n-1-j} x_3^{n-1-k} x_4^{\ell-2}  \right) \\
\\ = & \; \sigma_{n}  
\Delta(x_1,x_2,x_3,x_4)^{-1} \\ & \hspace{20.0pt} \times
\det\left( \begin{array}{cccc} 
\scs  (x_1-1)^3 \alpha_n(x_1) & \scs(x_2-1)^3\alpha_n(x_2) & \scs(x_3-1)^3\alpha_n(x_3) & \scs(x_4-1)^3 \alpha_n(x_4) \\
\scs  (x_1-1)^2 \beta_n(x_1) & \scs(x_2-1)^2 \beta_n(x_2) & \scs(x_3-1)^2\beta_n(x_3) & \scs(x_4-1)^2 \beta_n(x_4) \\
\scs  (x_1-1) \gamma_n(x_1) & \scs(x_2-1) \gamma_n(x_2) & \scs(x_3-1)\gamma_n(x_3) & \scs(x_4-1) \gamma_n(x_4) \\
\scs  \delta_n(x_1) & \scs\delta_n(x_2) & \scs\delta_n(x_3) & \scs\delta_n(x_4)
\end{array}\right)
\\
&- \sum_{t=1}^{4} \prod_{k=0}^{2} (1-x_{t+k}+x_{t+k-1} x_{t+k}) \;
 x_{t-1}^{n-2} x_{t+1}^{n-3} \;
 \mathcal A_{n-1}^{\textrm{TL}} \left(x_{t},\frac 1 
 {x_{t+1}} \right) \\
&-A_{n-2} \;(1-x_3+x_2x_3)(1-x_1+x_4 x_1)  \; x_2^{n-2}\; x_4^{n-2}\\
&-A_{n-2} \; (1-x_2+x_1x_2)(1-x_4+x_3x_4) \; x_1^{n-2}\; x_3^{n-2},
\end{split}
\ee
with the convention that $x_t = x_{t-4}$ for $5\le t\le 8$. The constant $\sigma_{n}$ is given by
\be \label{defsigman}
\sigma_n = \frac{1}{64 A_n^3 (2n-2)^{2}(2n-3)} 
\left(\frac{(n-2)!(3n-2)!}{(2n-3)!(2n-1)!}\right)^3.
\ee
\end{thm}

The basis for our enumeration results is a new expression for the partition function of the square ice model with the so-called ``domain wall'' boundary conditions at a special value of the parameter known as the ``crossing parameter''; see Theorem~\ref{thm:wrons} below.
This formula is of independent interest, both for its intrinsic theoretical value and because of its applicability to the problem of refined enumeration. As we learned from Filippo Colomo following the release of the preprint version of this paper, the formula is a close cousin of another expansion for the partition function derived by Colomo and Pronko \cite{colomopronko-orthogonal, colomopronko-emptiness}. A particular consequence of this relationship, brought to our attention by Colomo, is the following simplification of the above results.

\begin{thm}
\label{thm:tripquad-simplified}
The statements of Theorems~\ref{thm:trip} and \ref{thm:quad} above remain valid if the functions $\gamma_n(t)$ and $\delta_n(t)$ are replaced, in a manner analogous to the definition $\beta_n(t) = t\alpha_{n-1}(t)$, by
\be 
\begin{split} \label{redefalph}
\tilde{\gamma}_n(t) &= \mu_n t^2 \alpha_{n-2}(t), \\
\tilde{\delta}_n(t) &= \nu_n t^3 \alpha_{n-3}(t),
\end{split}
\ee
where $\mu_n$ and $\nu_n$ are suitable constants.
\end{thm}

Note that these functions are not identical to $\gamma_n(t)$ and $\delta_n(t)$, but the point is that the determinantal expansions \eqref{tripref} and \eqref{quadref} are unchanged by this modification.

In the Appendix we explain how to derive these simplified versions of the formulas \eqref{tripref}, \eqref{quadref} 
and discuss the relationship between our new formula for the square ice partition function (Theorem~\ref{thm:wrons}) and the Colomo-Pronko formula.

\bigskip
Note that the identities \eqref{tripref} and \eqref{quadref} (using either the original functions $\gamma_n, \delta_n$ or their simplified versions \eqref{redefalph} proposed by Colomo) make it possible to compute the coefficients $A_n^{\textrm{TLB}}(i,j,k)$ and $A_n^{\textrm{TLBR}}(i,j,k,\ell)$ in a computationally efficient manner using standard polynomial algebra. The Maple package \texttt{RefinedASM1234} accompanying this paper, which may be downloaded from the authors' web pages, provides a demonstration of the practical application of these formulas. 

It should also be noted that by making use of the generating function identities \eqref{tripref}, \eqref{quadref}, the coefficients $A_n^{TLB}(i,j,k)$ and $A_n^{TLBR}(i,j,k,\ell)$ can in principle be represented by explicit summation identities, similar to \eqref{eq:stroganov}. This is done by dividing out the factors outside the power series expansion on the left hand sides of \eqref{tripref} and \eqref{quadref}. The inverses of the factors, which are transferred to the right-hand side, are then expanded as double power series (using the expansion of $1/(1-x(1-y))$, which is easy to write down). The problem with this approach is that it will yield extremely complicated formulas: for example, in the case of the quadruply refined enumeration the formula for $A_n^{TLBR}(i,j,k,\ell)$, the formula will involve an 8-fold summation of binomial coefficients, unless some further ``magical'' simplification occurs. We conclude that the generating function identities \eqref{tripref} and \eqref{quadref} are in all likelihood the simplest ways of encoding the available information about the enumeration coefficients.

To conclude, we also note that, shortly after the release of the preprint version of our paper, a paper by Roger Behrend \cite{behrend} appeared on the arXiv repository in which he derives similar results to our main results. His paper in addition contains more general results concerning a further refinement of the quadruply refined enumeration based on two additional parameters, namely the inversion number of the ASM and the number of $-1$'s. See Section~3.3 of his paper for a detailed discussion of the relation of our results to his as well as an useful survey of the literature related to enumeration of ASMs and related combinatorial objects.

\paragraph{Acknowledgements.} We are grateful to Filippo Colomo for calling our attention to the connection of our results with the results from the papers \cite{colomopronko-orthogonal, colomopronko-emptiness}, and for suggesting the definitions \eqref{redefalph} that simplify our results. We also thank the anonymous referee for a careful reading of the paper and many insightful comments and suggestions.

\section{A new expression for the square ice partition function}
\label{sec:wrons}

The square ice model, also known as the six-vertex model, is a well known ``exactly solvable model'' from statistical physics \cite{baxter}.
The connection to ASMs appears when the model is considered with specific boundary conditions. We start with an $n \times n$ square lattice with $4n$ half-edges emanating from the boundary vertices. To the $i$th row, we associate a real-valued parameter $x_{i}$ and similarly, to the $j$th column, we associate a parameter~$y_{j}$. These are collectively called the \demph{spectral parameters}, with the $x_i$'s referred to as the \demph{row parameters} and the $y_j$'s termed the \demph{column parameters}.  A \demph{state}, or \demph{configuration}, of the model is an assignment of arrows to each edge (including the half-edges at the boundary)
so that every vertex has two incoming and two outgoing edges. There are six types of vertices, depicted in Figure~\ref{fig:sixvertex} below.

\begin{figure}[h]
\begin{center}
\begin{tabular}{cccccc}
\scalebox{0.19}{\includegraphics{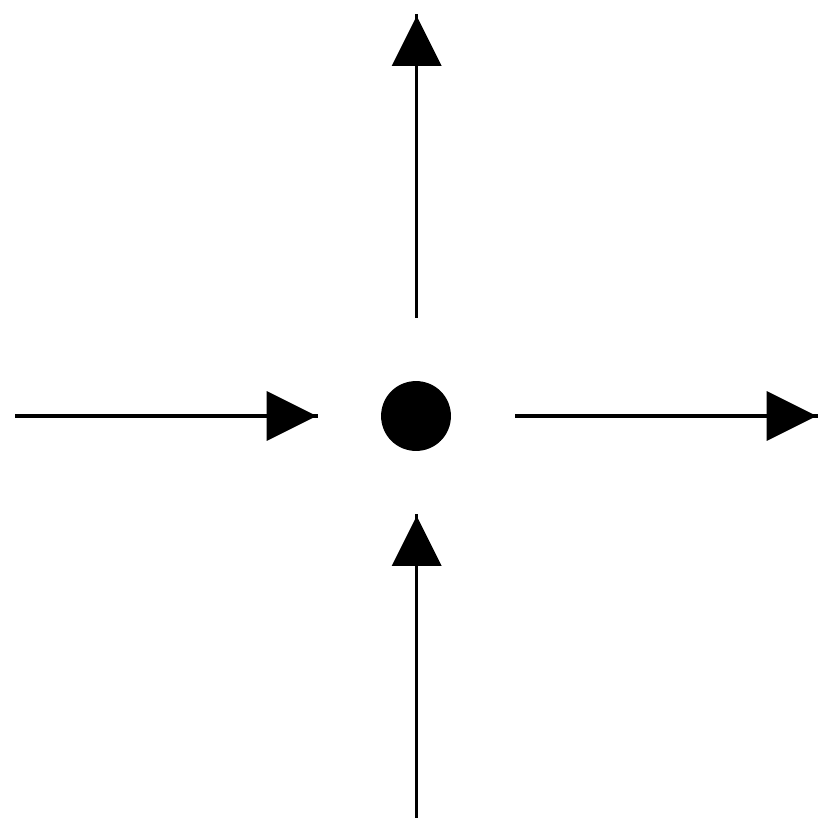}} &
\scalebox{0.19}{\includegraphics{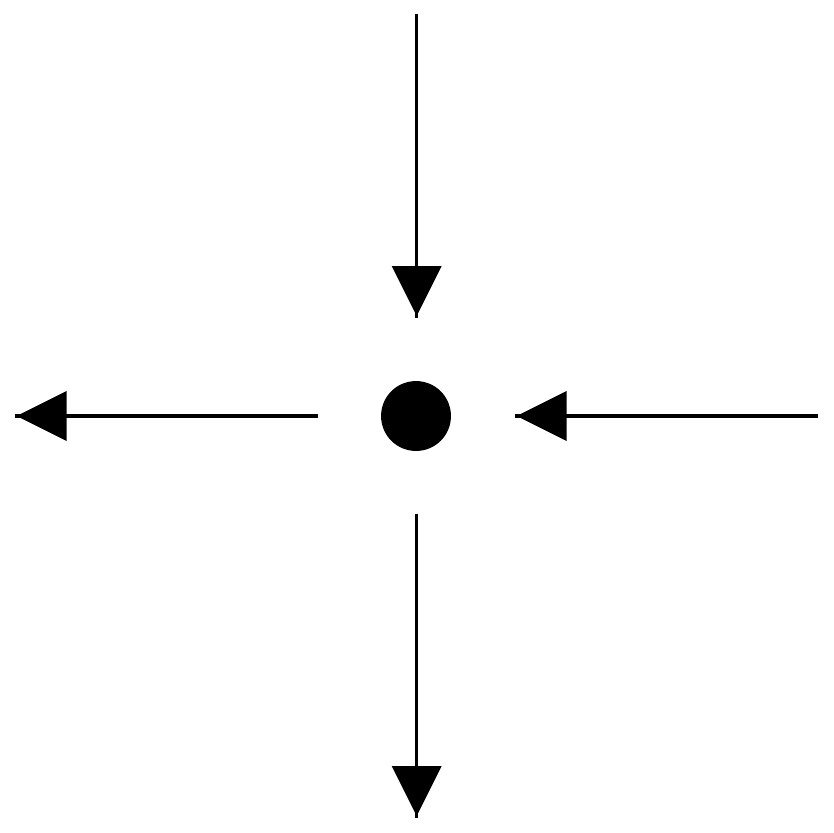}} &
\scalebox{0.19}{\includegraphics{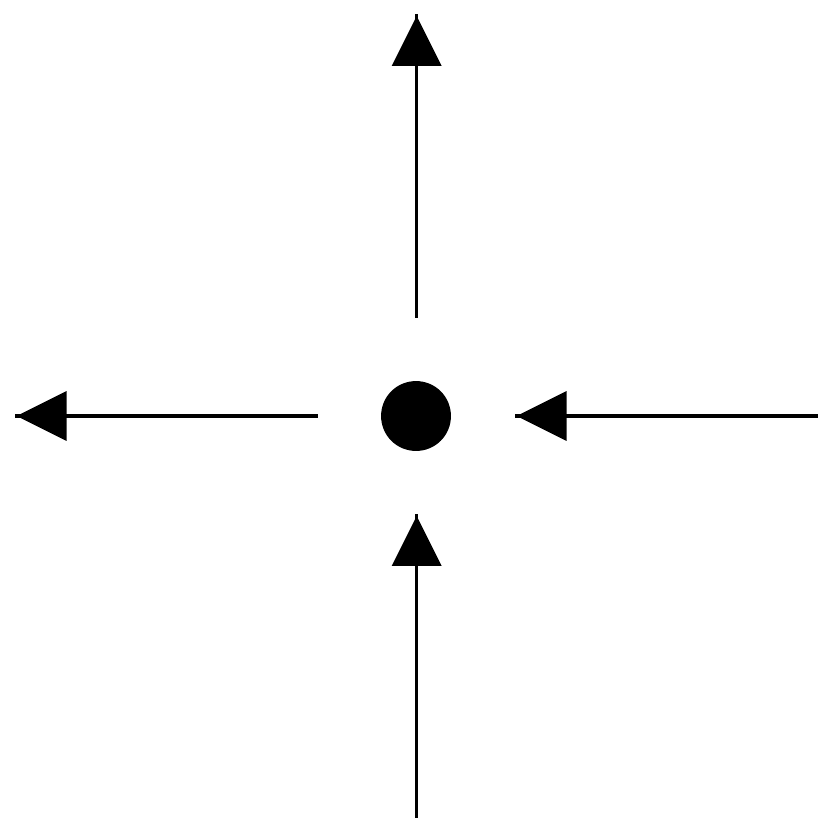}} &
\scalebox{0.19}{\includegraphics{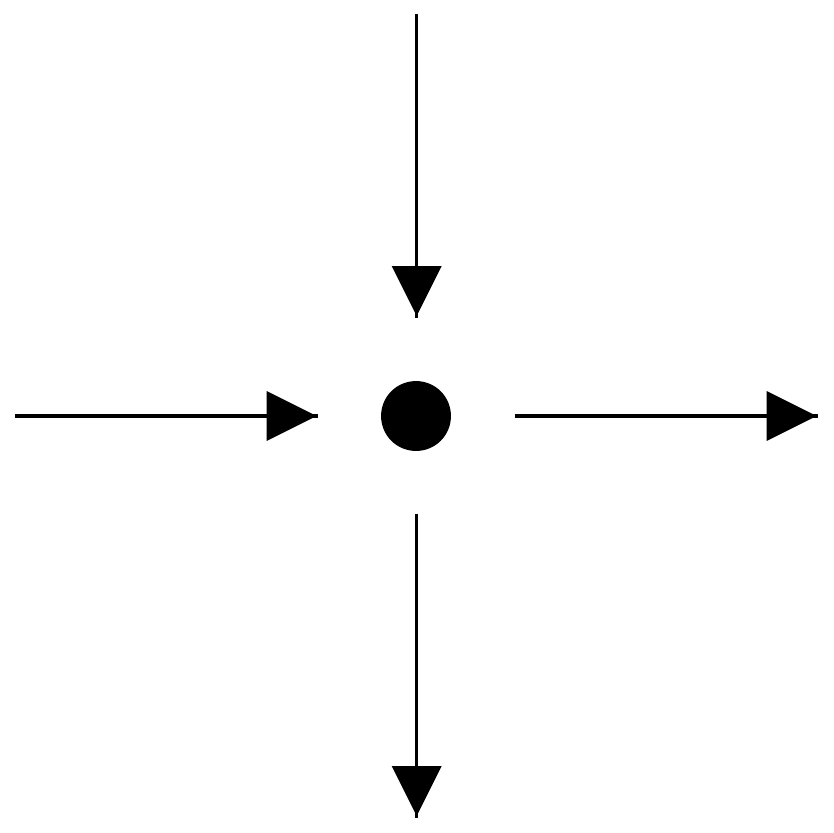}} &
\scalebox{0.19}{\includegraphics{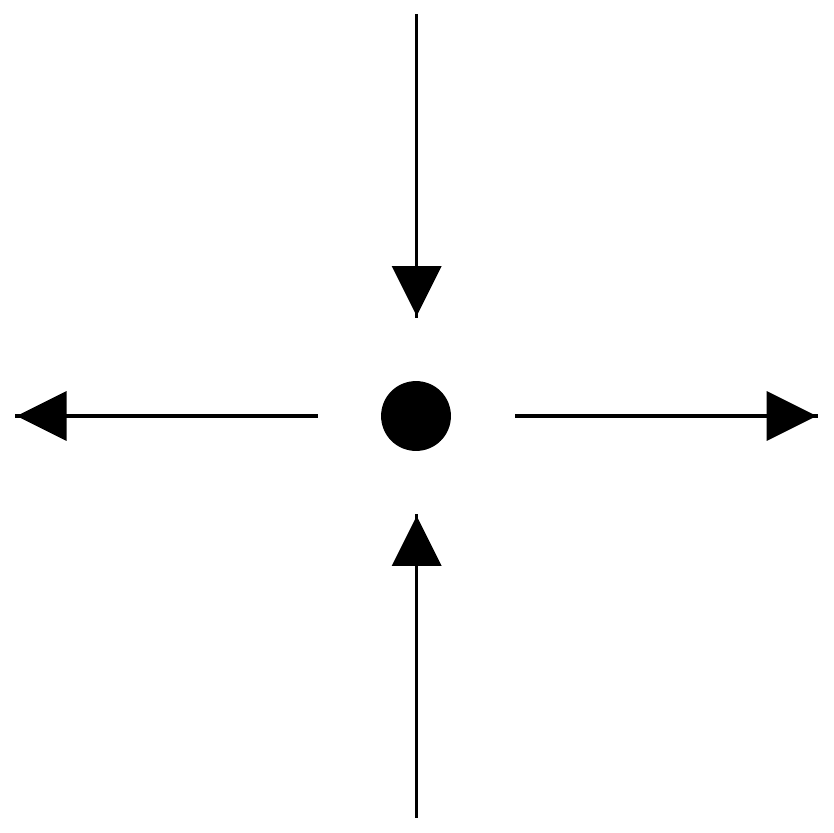}} &
\scalebox{0.19}{\includegraphics{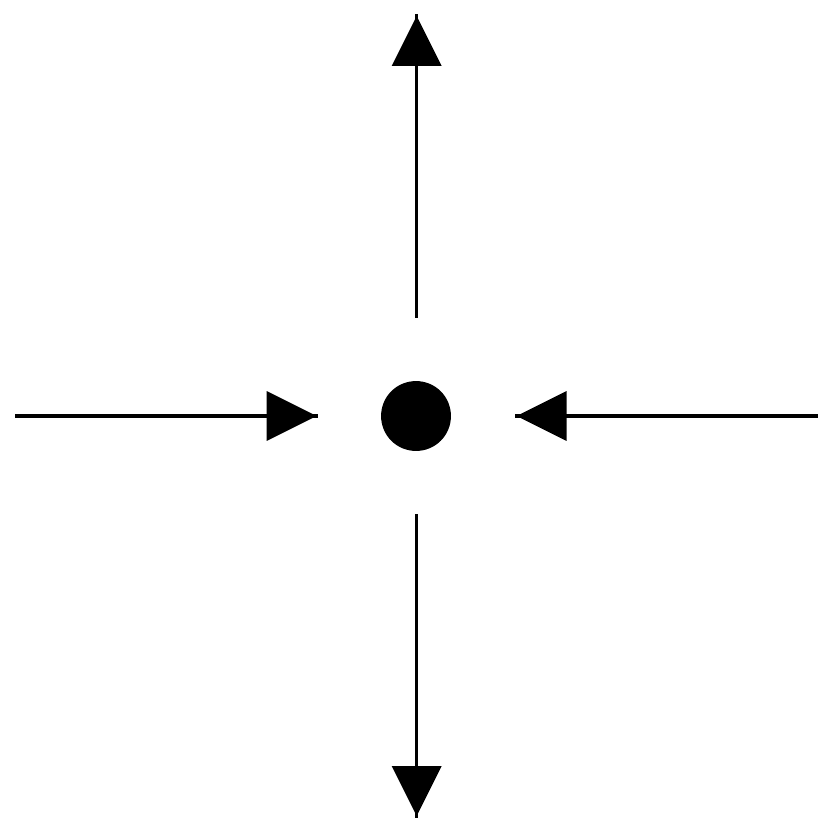}} \\
(type 1) &
(type 2) &
(type 3) &
(type 4) &
(type 5) &
(type 6)
\end{tabular}
\caption{The six types of vertices in a square ice configuration.
\label{fig:sixvertex}}
\end{center}
\end{figure}

\begin{figure}
\begin{center}
\scalebox{0.75}{\includegraphics{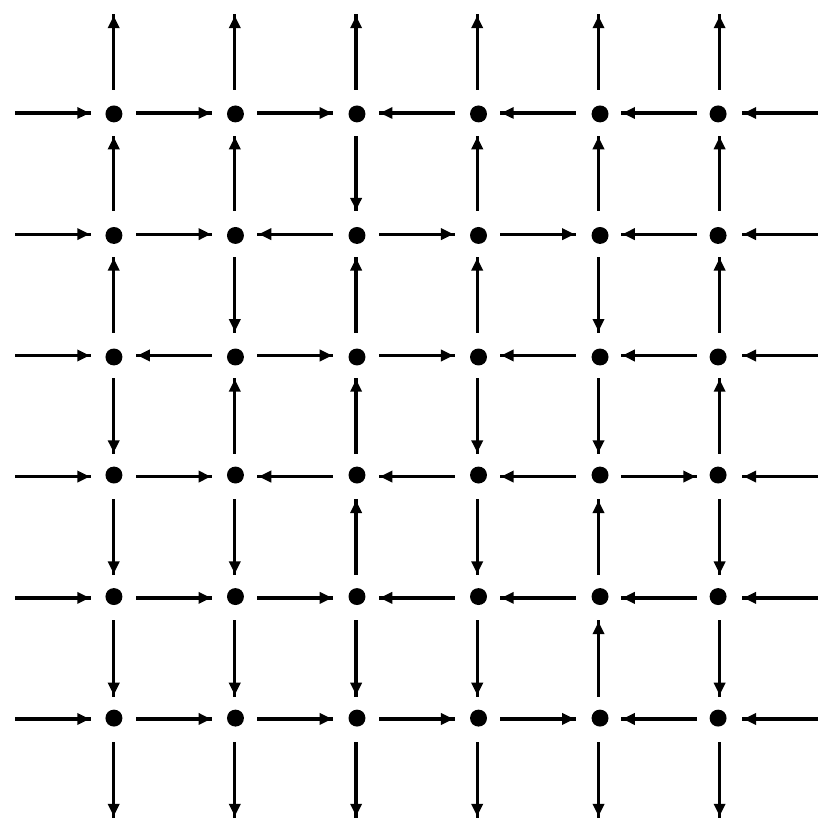}}
\caption{An example of a square ice configuration with domain wall boundary conditions\label{fig:sqiceeg}.}
\end{center}
\end{figure}

Type 2 is the rotation by $180^{\circ}$ of type 1, 
type 4 is the rotation by $180^{\circ}$ of type 3, and 
type 3 is the vertical reflection of type 1. Types 5 and 6 are rather special because
in both these cases, the two vertical (as also horizontal arrows) both point inwards or both point outwards. 
 The \demph{domain wall boundary conditions (DWBC)} correspond to the horizontal half-edges at the boundary pointing inwards and the vertical ones pointing outwards. See Figure~\ref{fig:sqiceeg} for an example. 
 
The relevance of square ice to the study of ASMs comes from the well-known fact that configurations of the square ice model with domain wall boundary conditions are in bijection with alternating sign matrices. Given a square ice configuration satisfying the DWBC, it can be translated to an alternating sign matrix by replacing every vertex of types 1 through 4 with a zero, vertices of type~5 with~$-1$'s and vertices of type 6 with 1's. As an example,
 the configuration in Figure~\ref{fig:sqiceeg} maps to the ASM in Figure~\ref{fig-asm-example} under this bijection.

Fix six functions $\mathfrak X_t(x,y)$, $1\le t\le 6$. 
Given a square ice configuration $C$ of order $n$, for each $1\le i,j\le n$ we assign a weight $\mathfrak X_{t_{i,j}}(x_{i},y_{j})$ to the vertex $(i,j)$, where $t_{i,j}$ is the type of the vertex at $(i,j)$. We then define the
weight of the configuration $C$ to be the product of the weights of all
the vertices,
\be
w(C) = \ds \prod_{1 \leq i,j \le n} \mathfrak X_{t_{i,j}}(x_{i},y_{j}).
\ee
Letting $\mathcal C_{n}$ denote the set of square ice configurations of order $n$, we define the \demph{partition function} as the sum of the weights over all elements of $\mathcal C_{n}$,
\be \label{defpartfn}
Z_{n} = \sum_{C \in \mathcal C_{n}} w(C),
\ee
which depends on $n$, the spectral parameters as well as the functions $\mathfrak X_{t}(\cdot,\cdot)$. In general, computing $Z_{n}$ is a hard problem, but something miraculous occurs when we set
\be \label{defx}
\begin{split}
\mathfrak X_{1}(x,y) &= \mathfrak X_{2}(x,y) 
= \frac{ \sin \left(\ds \tfrac \eta2+x-y \right)} {\ds \sin \eta}, \\
\mathfrak X_{3}(x,y) &= \mathfrak X_{4}(x,y) 
= \frac{ \sin \left(\ds \tfrac \eta2-x+y \right)} {\ds \sin \eta}, \\
\mathfrak X_{5}(x,y) &= \mathfrak X_{6}(x,y) = 1,
\end{split}
\ee
where $\eta$ is a real-valued parameter, known as the \demph{crossing parameter}.
In this case, the  partition function $Z_n=Z_{n}^{\eta}(\vec x; \vec y)$ can be expressed in terms of the  so-called \demph{Izergin-Korepin determinant} \cite{izercokkor}, namely as
\be \label{izergin-korepin}
\begin{split}
Z_{n}^{\eta}(\vec x; \vec y) &= (-1)^{\binom n2} \frac 1{\sin^{n(n-1)} (\eta) }
\frac{ \ds
\prod_{1 \leq i,j \leq n} \sin \left(x_{i}-y_{j} + \tfrac \eta2 \right)
\sin \left(x_{i}-y_{j} - \tfrac \eta2 \right)}
{\ds \prod_{1 \leq i < j \leq n} \sin \left(x_{i}-x_{j}\right)
\sin \left(y_{i}-y_{j} \right)} 
\\ & \ \ \ \times 
\Det{i,j=1}{n} M(\vec x; \vec y),
\end{split}
\ee
where $M(\vec x; \vec y)$ is the $n \times n$ square matrix with entries given by
\be
M(\vec x; \vec y)_{i,j} = \frac 1{\sin \left(x_{i}-y_{j} + \tfrac \eta2 \right)
\sin \left(x_{i}-y_{j} - \tfrac \eta2 \right)},
\ \quad (1\le i,j\le n).
\ee

The partition function is clearly a trigonometric polynomial, which, by \eqref{izergin-korepin}, is seen to be a symmetric function of the $x_{i}$'s and (independently) of the $y_{i}$'s. Moreover, it is not difficult to see from \eqref{defpartfn} that it is of degree at most $n-1$ in each variable. A rather surprising fact about the partition function, shown by Stroganov \cite{stroganov}, is that when the crossing parameter $\eta$ is set at $\eta=2\pi/3$, the partition function becomes completely symmetric in the variables $\{\cup x_{i} \bigcup  \cup y_{i} \}$. 
Stroganov used this insight to find yet another determinental expression (essentially a Schur function, see \cite[Section 2.5.6]{zinnjustin}) for the
partition function, which led to an alternate proof of the Refined ASM Theorem and to the formulas 
\eqref{eq:stroganov}, \eqref{eq:stroganov2} for the ``top-bottom'' and ``top-left'' doubly refined enumerations.

From now on, we assume that the weights $\mathfrak{X}_t$ are given by \eqref{defx} and that $\eta=2\pi/3$ (this point in the parameter space of the square ice model is sometimes referred to as the \demph{combinatorial point}). Since the partition function is completely symmetric in the $x$~and~$y$ variables, following Stroganov we use a uniform notation for the spectral parameters, viz.\ $u_{i}$ for $i \in \{1,\dots,2n\}$. To avoid excessive clutter of notation, we reuse $Z_{n} = Z_{n}(\vec u)$ to denote the partition function used in this way.
We also use the notation $Z_{n}(u_{1},\dots,u_{k}, 0^{2n-k})$ for the partition function with the last $2n-k$ variables set to zero. 

Let 
\be \label{defabc}
\begin{split}
\weighta(u) &= \frac{2}{\sqrt{3}} \sin\left(\frac\pi3 + u \right), \\
\weightb(u) &= \frac{2}{\sqrt{3}} \sin\left(\frac\pi3 - u \right), \\
c(u) &= 1,
\end{split}
\ee
and note that the $\mathfrak{X}_t$ may be expressed in terms of the functions $a(\cdot), b(\cdot)$ and $c(\cdot)$ as
\be
\begin{split}
\mathfrak{X}_1(x,y) &= \mathfrak{X}_2(x,y) = \weighta(x-y), \\
\mathfrak{X}_3(x,y) &= \mathfrak{X}_4(x,y) = \weightb(x-y), \\
\mathfrak{X}_5(x,y) &= \mathfrak{X}_6(x,y) = 1 = c(x-y). \\
\end{split}
\ee
Note that $a(0)=b(0)=c(0)=1$, so in particular when all the spectral parameters are set to $0$, the weight of every square ice configuration is equal to $1$, and we get that $Z_n(0^{2n}) = A_n$, the total number of ASMs of order~$n$.

The following identities are easily proved and will be of use to us below:
\be \label{abidentities}
\begin{split}
\sin(u) &= \frac{\sqrt{3}}{2}(\weighta(u)-\weightb(u)), \\
\cos(u) &= \frac12 (\weighta(u)+\weightb(u)), \\
\cos\left(\frac\pi3 + u \right) &= \frac12 (-\weighta(u)+2\weightb(u)), \\
\cos\left(\frac\pi3 - u \right) &= \frac12 (2\weighta(u)-\weightb(u)), \\
\sin(u-v) &= \frac{\sqrt{3}}{2}(\weighta(u) \weightb(v) - \weightb(u) \weighta(v)), \\
\sin\left(\frac\pi3+u-v \right) &= \frac{\sqrt{3}}{2}(\weighta(u) \weighta(v)+\weightb(u) \weightb(v) - \weightb(u) \weighta(v)).
\end{split}
\ee
Now, following \cite{stroganov}, let
\be \label{deffn}
f_{n}(u) = \sin^{2n-1}(u) Z_{n}(u, 0^{2n-1}),
\ee
and note that
\begin{align}
f_n(0)& = f_n'(0) = f_n''(0) = \ldots = f_n^{(2n-2)}(0) = 0, \label{fn-derivatives1} \\
f_n^{(2n-1)}(0) &= (2n-1)! Z_n(0^{2n}) = (2n-1)! A_n. \label{fn-derivatives2}
\end{align}
The function $f_n(u)$ encodes information about the usual ``top row'' refined enumeration of ASMs: by a standard ``first row expansion'' (dividing square ice configurations into classes according to their first row, and noting that in the specialization $Z_n(u,0^{2n-1})$ the contribution to the weight from all other rows is constant), we see that it may be written as
\be \label{firstrowexpansion}
f_n(u) = \sin^{2n-1}(u)  \sum_{k=1}^{n} A_{n,k} \weighta^{k-1}(u) \weightb^{n-k}(u)
\ee
(see \cite[Section 4]{stroganov}). Stroganov also found a fairly simple explicit expansion for $f_n$ as a trigonometric polynomial:

\begin{lem} \label{lem:kappan}
\be \label{stroganovexpansion}
f_n(u) =
\kappa_n
\sum_{m=0}^{n-1} \binom{n-\tfrac43}{m}\binom{n-\tfrac23}{n-m-1} \sin\big( (4-3n+6m)u \big),
\ee
where the constant $\kappa_n$ is given by
\be
\label{kappan-value}
\kappa_n = \left(\frac{-3}{4}\right)^{n-1} \frac{A_{n-1}}{\binom{2n-2}{n-1}}.
\ee
\end{lem}

\begin{proof} The formula \eqref{stroganovexpansion} was proved by Stroganov (who used it to obtain a new derivation of the formula for $A_{n,k}$), except for the value of the constant $\kappa_n$. To compute $\kappa_n$, we set $u=\pi/3$ and evaluate $f_n(\pi/3)$ in two ways. On the one hand, from \eqref{deffn} we have that
\begin{align*}
f_n(\pi/3) &= \sin^{2n-1}\left(\pi/3\right) \sum_{k=1}^n A_{n,k} a(\pi/3)^{k-1} b(\pi/3)^{n-k} \\ &=
\left(\frac{\sqrt{3}}{2}\right)^{2n-1} A_{n,n} = \frac{2}{\sqrt{3}} \left(\frac34\right)^n A_{n-1}
\end{align*}
(since $a(\pi/3)=1, b(\pi/3)=0$); on the other hand, \eqref{stroganovexpansion} gives
\begin{align*}
f_n(\pi/3) &= \kappa_n \sum_{m=0}^{n-1} \binom{n-\tfrac43}{m} \binom{n-\tfrac23}{n-m-1} \sin\left(-\pi n+\tfrac{4\pi}{3}\right) \\ &= (-1)^{n-1} \frac{\sqrt{3}}{2} \kappa_n \sum_{m=0}^{n-1}
\binom{n-\tfrac43}{m} \binom{n-\tfrac23}{n-m-1} \\ &=
(-1)^{n-1} \frac{\sqrt{3}}{2} \kappa_n \binom{2n-2}{n-1},
\end{align*}
where the last step follows from the Chu-Vandermonde summation identity. Comparing the two expressions gives \eqref{kappan-value}.
\end{proof}

As a side remark, note that one can get a different expression for $\kappa_n$ by again equating the two expressions \eqref{deffn} and \eqref{stroganovexpansion} and taking the limit as $u\to 0$
(using L'H\^opital's rule). This gives
$$ \kappa_n = (-1)^{n-1} (2n-1)! A_n \left(\sum_{m=0}^{n-1} \binom{n-\tfrac43}{m}\binom{n-\tfrac23}{n-m-1} (4-3n+6m)^{2n-1}\right)^{-1}.
$$
The fact that this is equal to the right-hand side of \eqref{kappan-value} implies the following interesting identity.

\begin{cor}
The summation identity
\be \label{summation-identity}
\sum_{m=0}^{n} \binom{n-\tfrac13}{m}\binom{n+\tfrac13}{n-m} (1-3n+6m)^{2n+1} = \left( \frac 43 \right)^n \frac{(3n+1)!}{n!}.
\ee
holds for all $n\ge 1$.
\end{cor}

An independent proof of \eqref{summation-identity} was found by Amdeberhan \cite{amdeberhan}.

Our next result will give a hint that the remarkable function $f_n(u)$ also holds the key to higher order refined enumerations of ASMs, which are related to the $k$-variate specialization $Z_n(u_1,\ldots,u_k, 0^{2n-k})$.

\begin{thm}[Partition function formulas] \label{thm:wrons}
The partition function with all but $k$ variables set to zero can be written as
\be \label{wrons}
Z_{n}(u_{1},\dots,u_{k}, 0^{2n-k}) = \zeta_{n,k} 
\frac { \ds \Det{i,j=1}{k} \left( \frac{d^{j-1} f_{n}(u_{i})}{du_{i}^{j-1}} \right)}
{\ds \prod_{i=1}^{k} \sin^{2n-k}(u_{i})
\prod_{1 \leq i < j \leq k} \sin(u_{i}-u_{j})},
\ee
where the constant $\zeta_{n,k}$ is given by
\be \label{eq:formula-zetank}
\zeta_{n,k}= \frac 1{A_{n}^{k-1} \prod_{j=1}^{k-1}(2n-j)^{k-j}}.
\ee
In particular, taking $k=2n$ we get the expression
\be
Z_n(u_1,\ldots,u_{2n}) = \frac 1{A_{n}^{2n-1} H(2n-1)}
\ds \Det{i,j=1}{2n}\! \left( \frac{d^{j-1} f_{n}(u_{i})}{du_{i}^{j-1}} \right)
\!\!\prod_{1 \leq i < j \leq 2n}\!\!\! \sin(u_{i}-u_{j})^{-1}
\ee
for the full partition function, where $H(n)=\prod_{k=1}^n k^k$ is the hyperfactorial.
\end{thm}

We note in passing that several different representations for the partition function have been found, see, e.g., \cite{izercokkor,colomopronkopartfn1,
colomopronkopartfn2,galleas, behrend-pdf-pzj}.

The proof will require the following lemma.

\begin{lem} \label{lem-systemeq}
Let $r_1,\ldots,r_N$ be distinct real numbers, and let $0\le k< N$. The solutions to the homogeneous system of linear equations
\begin{align}
\sum_{j_1=1}^N z_{j_1,\ldots,j_k} r_{j_1}^m &= 0  \qquad (1\le j_2,\ldots,j_k\le N,\ 0\le m\le N-k-1), \label{eq:system-rowmoments} \\
z_{\sigma(j_1),\ldots,\sigma(j_k)} &= \operatorname{sgn}(\sigma)z_{j_1,\ldots,j_k} \qquad  (\sigma\in S_k),
\label{eq:system-antisym}
\end{align}
in the family of unknowns $(z_{j_1,\ldots,j_k})_{1\le j_1,\ldots,j_k\le N}$, are determined up to a global multiplicative constant.
\end{lem}

\begin{proof}
The anti-symmetry condition \eqref{eq:system-antisym} implies that any variable $z_{j_1,\ldots,j_k}$ where $j_1,\ldots,j_k$ are not all distinct is $0$. 
It will therefore be enough to show that for any two index vectors $(j_1,\ldots,j_k), (j_1',\ldots,j_k')$ where $j_1,\ldots,j_k$ are distinct and $j_1',\ldots,j_k'$ are distinct, the ratio $z_{j_1,\ldots,j_k}/z_{j_1',\ldots,j_k'}$ is determined uniquely by the equations. Furthermore, if we prove this claim in the case that the vectors $(j_1,\ldots,j_k), (j_1',\ldots,j_k')$ differ in only one position then the more general claim follows by successively replacing one index at a time (using anti-symmetry to permute the order of the indices). Again because of anti-symmetry, it suffices to show this when the vectors differ only in the \emph{first} position. That is, we wish to show that the ratio $z_{j_1,\ldots,j_k}/z_{j_1',j_2,\ldots,j_k}$ is uniquely determined where $j_1,j_1',j_2,\ldots,j_k$ are distinct. Fixing distinct indices $j_2,\ldots,j_k$, we can consider the equations \eqref{eq:system-rowmoments} as $m$ ranges from $0$ to $N-k-1$ as a system of $N-k$ linear equations with $N-k+1$ unknowns (think of the sum on the left-hand side of \eqref{eq:system-rowmoments} as ranging only over $j_1\in \{1,\ldots,N\}\setminus\{j_2,\ldots,j_k\}$, since the $z$-coefficients corresponding to the skipped values are already known to vanish) . The coefficient matrix of this system is $(r_j^m)_{0\le m\le N-k-1,\ j\in\{1,\ldots,N\}\setminus\{j_2,\ldots,j_k\}}$. This matrix has rank $N-k$, since the $(N-k)\times(N-k)$ minor which consists of the first $N-k$ columns (or indeed any $N-k$ columns) is a Vandermonde determinant, which is non-zero since the $r_i$ are distinct. It follows that the kernel of the matrix is one dimensional, which implies our claim.
\end{proof}

\begin{proof}[Proof of Theorem~\ref{thm:wrons}]
Define the function $g_{n,k}$ as
\be
\begin{split}
g_{n,k}(u_{1},\dots,u_{k}) &= Z_{n}(u_{1},\dots,u_{k}, 0^{2n-k})
\\ & \ \ \ \ \times 
\prod_{i=1}^{k} \sin^{2n-k}(u_{i})
\prod_{1 \leq i < j \leq k} \sin(u_{i}-u_{j}).
\end{split}
\ee
The proof is based on the observation that $g_{n,k}$ satisfies
the following list of properties: 
\begin{enumerate}

\item $g_{n,k}$ is a trigonometric polynomial in the $u_{i}$'s of degree at most $3n-2$
in each variable.

\item $g_{n,k}(\dots,u_{i}+\pi,\dots) = (-1)^{n} g_{n,k}(\dots,u_{i},\dots)$.

\item $g_{n,k}(\dots,u_{i},\dots) + g_{n,k}(\dots,u_{i}+\frac{2\pi}3,\dots)
+ g_{n,k}(\dots,u_{i}+\frac{4\pi}3,\dots) = 0$.

\item $g_{n,k}$ is divisible by $ \sin^{2n-k}(u_{i})$ for each $i$.

\item $g_{n,k}(\dots,u_{i},\dots,u_{j},\dots) = -g_{n,k}(\dots,u_{j},\dots,u_{i},\dots)$.

\end{enumerate}

Property 1 follows from the fact that the partition function $Z_n$ is a trigonometric polynomial of degree $\le n-1$ in each variable. Property 2 follows from the relation
\be \label{eq:zn-parity}
Z_n(\dots,u_{i}+\pi,\dots) = (-1)^{n-1} Z_n(\dots,u_{i},\dots)
\end{equation}
satisfied by $Z_n$, which is a consequence of the fact that $a(u+\pi)=-a(u)$, $b(u+\pi)=-b(u)$, together with the observation that each row and column of a square ice configuration contain an odd number of vertices with weight $c$ (translate this to the language of ASMs to see why).
Property 3 is a consequence of the symmetry of $Z_n$ in the variables $u_1,\ldots,u_{2n}$ together with the following functional equation proved by Stroganov \cite{stroganov}:
\be \label{eq:stroganov-functionaleq}
F_n(u_1,\ldots)+F_n\left(u_1+\tfrac{2\pi}{3},\ldots\right)
+F_n\left(u_1+\tfrac{4\pi}{3},\ldots\right)=0,
\ee
where $F_n(u_1,\ldots,u_{2n})$ is defined by
$$
F_n(u_1,\ldots,u_{2n}) = Z_n(u_1,\ldots,u_{2n}) \prod_{j=2}^{2n} \sin(u_1-u_j).
$$
Properties 4 and 5 are immediate from the definition of $g_{n,k}$.

We will now show that any function which satisfies these five properties is determined uniquely up to a normalization constant depending only on $n$ and $k$. The idea is a generalization of Stroganov's approach to the proof of the Refined ASM theorem \cite{stroganov}.
Starting from the first property, we express $g_{n,k}$ in canonical form,
\be \label{eq:gnk-fourier}
g_{n,k}(u_{1},\dots,u_{k}) = \sum_{a_{1}=-3n+2}^{3n-2} \ldots
 \sum_{a_{k}=-3n+2}^{3n-2} c_{a_{1},\dots,a_{k}} \exp \left(\sqrt{-1} \sum_{j=1}^{k} a_{j} u_{j}\right),
\ee
where the Fourier coefficients $c_{a_{1},\dots,a_{k}}$ are to be determined. 

Using the second property, we see that $c_{a_{1},\dots,a_{k}}=0$ whenever some $a_{j}$ does not have the same parity as $n$. Using the third property, $c_{a_{1},\dots,a_{k}}=0$ whenever some
$a_{i}$ is divisible by $3$. So we may rewrite \eqref{eq:gnk-fourier} as
$$ 
g_{n,k}(u_{1},\dots,u_{k}) = \sum_{a_{1}\in I_n} \ldots
 \sum_{a_{k}\in I_n} c_{a_{1},\dots,a_{k}} \exp \left(\sqrt{-1} \sum_{j=1}^{k} a_{j} u_{j}\right),
$$
where $$I_n = \{ 2-3n \le a\le 3n-2\,:\, a\equiv n\ (\operatorname{mod}2), a\equiv 1,2\ (\operatorname{mod}3) \}.$$
Note that $|I_n| = 2n$. Now, the fourth property implies that
$$ \frac{\partial ^m}{\partial u_1^m} \bigg|_{u_1=0} g_{n,k}(u_1,\ldots,u_k) = 0, \qquad (0\le m\le 2n-k-1), $$
which translates to the equations
\be \label{eq:system-rowmoments2}
\sum_{a_1 \in I_n} c_{a_1,\ldots,a_k} a_1^m = 0
\qquad (a_2,\ldots,a_k\in I_n, \ 0\le m\le 2n-k-1)
\ee
on the coefficients. Similarly, the anti-symmetry property 5 implies the corresponding anti-symmetry condition
\be \label{eq:system-antisym2}
c_{\sigma(a_1),\ldots,\sigma(a_k)} = \operatorname{sgn}(\sigma)c_{a_1,\ldots,a_k}
\ee
on the coefficients, which holds for any permutation $\sigma$ of the elements of $I_n$. But we now see that \eqref{eq:system-rowmoments2} and \eqref{eq:system-antisym2} are identical to the system of equations \eqref{eq:system-rowmoments}, \eqref{eq:system-antisym} from Lemma~\ref{lem-systemeq}, with the substitution $N=2n$ and with $r_1,\ldots,r_N$ denoting the elements of $I_n$. The conclusion of Lemma~\ref{lem-systemeq} guarantees that the coefficients $c_{a_1,\ldots,a_k}$, and hence the function $g_{n,k}$, are determined uniquely up to a multiplicative constant, which was our claim up to this point.

Next, we claim that 
\be
h_{n,k}(u_{1},\dots,u_{k}) = \Det{i,j=1}{k}\left( \frac{d^{j-1} f_{n}(u_{i})}{du_{i}^{j-1}} \right)
\ee
satisfies the same list of properties discussed above. First of all, $f_{n}(u)$ is a trigonometric polynomial of degree $3n-2$ in $u$ and taking its derivatives cannot increase the degree, so $h_{n}$ satisfies the first property. The second and third properties follow immediately from the equations
\begin{align*}
f_n(u+\pi)&=(-1)^n f_n(u), \\
0&=f_n(u)+f_n\left(u+\frac{2\pi}{3}\right)+f_n\left(u+\frac{4\pi}{3}\right)
\end{align*}
satisfied by $f_n(u)$ (the first equation is again a consequence of \eqref{eq:zn-parity}, and the second equation is the special case $u_2=\ldots=u_{2n}=0$ of \eqref{eq:stroganov-functionaleq}).

The fourth property holds because the smallest factor of $\sin(u)$ in a derivative of $f_{n}$ is $\sin^{2n-k}(u)$ in the $(k-1)$th derivative.
The fifth and last property is immediate. Thus, we have shown that
$g_{n}(u_{1},\dots,u_{k}) = \zeta_{n,k} h_{n}(u_{1},\dots,u_{k})$ where $\zeta_{n,k}$ is a proportionality constant.

It remains to prove the formula \eqref{eq:formula-zetank} for $\zeta_{n,k}$. We shall do so by an induction on $k$ for each fixed $n$. We first note that $\zeta_{n,1}=1$ 
by \eqref{deffn}. We now evaluate $\zeta_{n,k}$ by taking the limit of all $u_{i}$'s to zero in \eqref{wrons},
\be
\frac{\zeta_{n,k}}{A_{n}} = \lim_{u_{1},\dots,u_{k} \to 0}
\frac {\ds \prod_{i=1}^{k} \sin^{2n-k}(u_{i})
\prod_{1 \leq i < j \leq k} \sin(u_{i}-u_{j})}
{ \ds \Det{i,j=1}{k} \left( \frac{d^{j-1} f_{n}(u_{i})}{du_{i}^{j-1}} \right)}.
\ee
We compute the limit as an iterative limit, by taking the variables successively to $0$; in fact, because of the inductive nature of the computation we only need to let $u_k\to 0$ since a recursive structure is revealed. To do this, we use L'H\^opital's rule by taking the $(2n-k)$th partial derivative of the numerator and denominator with respect to $u_k$ and setting each to $0$. Consider first the relevant terms in the numerator. We have to evaluate
\be
\frac{\partial ^{2n-k}}{\partial u_{k}^{2n-k}}\textrm{\raisebox{-5pt}{$\bigg|_{u_{k}=0}$}} \left( \sin^{2n-k}(u_{k}) 
\prod_{i=1}^{k-1} \sin(u_{i}-u_{k}) \right) .
\ee
This multiple derivative will lead to a sum involving multinomial coefficients over all possible combinations of derivatives of the $k$ terms inside in general. However, since we are going to set $u_{k}=0$ at the very end, the only term that contributes is the one where all $2n-k$ derivatives act on $\sin^{2n-k}(u_{k})$. We thus obtain, after easy manipulations,
\be
(2n-k)! \prod_{i=1}^{k-1} \sin(u_{i}).
\ee
We now consider the denominator. Since $u_{k}$ is only present in the last row of the matrix, we can perform a cofactor expansion of the determinant using the last row. 
By \eqref{fn-derivatives1},
the only term that contributes after taking the $(2n-k)$th derivative is the cofactor of the entry in position $(k,k)$, which leads to
\begin{align*} 
\Bigg(\frac{d^{2n-1}}{du_{k}^{2n-1}} &\textrm{\raisebox{-5pt}{$\bigg|_{u_{k}=0}$}}  f_{n}(u_{k}) \Bigg) \times
\Det{i,j=1}{k-1} \left( \frac{d^{j-1} f_{n}(u_{i})}{du_{i}^{j-1}} \right)
\\ &= (2n-1)!A_n \Det{i,j=1}{k-1} \left( \frac{d^{j-1} f_{n}(u_{i})}{du_{i}^{j-1}} \right)
\end{align*}
(using \eqref{fn-derivatives2}).
To summarize, we have shown that
\be
\begin{split}
\frac{\zeta_{n,k}}{A_{n}} 
&= \frac{(2n-k)!}{ (2n-1)! A_{n}}
 \lim_{u_{1},\dots,u_{k-1} \to 0} 
 \frac {\ds \prod_{i=1}^{k-1} \sin^{2n-k+1}(u_{i})
\prod_{1 \leq i < j \leq k-1} \sin(u_{i}-u_{j})}
{(2n-1)! A_{n}
\Det{i,j=1}{k-1} \left( \frac{d^{j-1} f_{n}(u_{i})}{du_{i}^{j-1}} \right)}, \\
&= \frac{(2n-k)!}{ (2n-1)! A_{n}} \frac{\zeta_{n,k-1}}{A_{n}}.
\end{split}
\ee
It is easy to check that the expression on the right-hand side of \eqref{eq:formula-zetank} satisfies this.
\end{proof}

When one is interested in refined enumeration formulas, one would like to compute 
$Z_{n}(u_{1},\dots,u_{k},0^{2n-k})$ for large $n$ and fixed $k$. In this case 
\eqref{wrons} is very useful because the order of the determinant is fixed. Note, however, that this formula would not be very practical to use if the
computation of the derivatives of $f_{n}$ were not efficient enough.
As we shall see in the next section, for the purpose of deriving enumeration formulas we are interested in a representation of the derivatives $f_n^{(m)}(u)$ not as a standard trigonometric polynomial (this would be immediate from \eqref{stroganovexpansion}), but instead we need a specific expansion involving products of powers of $a(u)$ and $b(u)$ factors, similarly to \eqref{firstrowexpansion}, with an additional factor consisting of a power of $\sin(u)$.
We now show that there is an efficient recursive algorithm for the computation of such a representation.

\begin{prop} \label{prop:fnmrecur}
For $0\le m\le 2n-1$, the $m$th derivative of $f_{n}$ can be written as
\be \label{fnmrecur}
f^{(m)}_{n}(u) = \sin^{2n-1-m}(u) \sum_{k=1}^{n+m} c_{n,m,k} a^{k-1}(u) b^{n-k+m}(u),
\ee
where $c_{n,0,k} = A_{n,k}$, and for $m>0$ we have the recurrence
\be \label{cnmkrecur}
\begin{split}
c_{n,m,k} = &\tfrac12 \Big( 
- 2kc_{n,m-1,k+1} 
 + (n+4k-2m-2) c_{n,m-1,k}  \\
&  + (5n-4k+2m+2) c_{n,m-1,k-1} 
-2(n-k+m+1) c_{n,m-1,k-2} 
\Big) .
\end{split}
\ee
\end{prop}

\begin{proof}

Since the argument of $a$ and $b$ will always be $u$, we will omit it in this proof. Starting with 
\be
f^{(m-1)}_{n}(u) = \sin^{2n-m}(u) \sum_{k} c_{n,m-1,k} a^{k-1} b^{n+m-1-k}
\ee
and differentiating, we get
\be
\begin{split}
f^{(m)}_{n}(u) =& (2n-m) \sin^{2n-1-m}(u) \cos(u)
 \sum_{k} c_{n,m-1,k} a^{k-1} b^{n+m-1-k} \\
&+  \sin^{2n-m}(u) \sum_{k} c_{n,m-1,k} 
\left[ (k-1) a^{k-2} b^{n+m-1-k} \frac{2}{\sqrt{3}}\cos\left(\frac\pi3 + u \right) \right. \\
&- \left. (n+m-1-k) a^{k-1} b^{n+m-2-k} \frac{2}{\sqrt{3}}\cos\left(\frac\pi3 - u \right) \right].
\end{split}
\ee
Now we take out the common factor of $\sin^{2n-1-m}(u)$ and use the identities \eqref{abidentities} to replace all trigonometric factors by sums of products of $a$'s and 
$b$'s. After simplifying and collecting terms, we get
\be
\begin{split}
f^{(m)}_{n}(u) =& \tfrac12 \sin^{2n-1-m}(u) 
 \sum_{k} c_{n,m-1,k}  \left[  -2(k-1) a^{k-2} b^{n-k+m+1}  \right. \\
&+ (n-2m-2+4k) a^{k-1} b^{n+m-k} \\
&+ (5n+2m-2-4k) a^{k} b^{n+m-1-k} \\
&- \left. 2(n+m-1-k) a^{k+1} b^{n+m-2-k} \right],
\end{split}
\ee
which is another way of encoding \eqref{cnmkrecur}. Note that if we assume inductively that the range of values of $k$ for which the coefficient $c_{n,m-1,k}$ is nonzero is $1\le k\le n+m-1$, then the range of $k$'s over which we should consider $c_{n,m,k}$ is $0\le k\le n+m+1$. However, the coefficient $2k$ appearing in front of the factor $c_{n,m-1,k+1}$ in \eqref{cnmkrecur} ensures that $c_{n,m,0}=0$, and similarly the coefficient $2(n-k+m+1)$ of $c_{n,m-1,k-2}$ ensures that $c_{n,m,n+m+1}=0$, so in fact we only get nonzero coefficients $c_{n,m,k}$ when $1\le k\le n+m$. This explains the range of summation in \eqref{fnmrecur}.
\end{proof}

The proof of Proposition~\ref{prop:fnmrecur} suggests a new way of looking at
the partition function $Z_{n}$. The usual way of considering it as a trigonometric polynomial has been replaced by a polynomial in $a$ and $b$. This is in some sense simpler because the total degree of the polynomial is constant. Underlying this is a change-of-basis transformation in the space of trigonometric polynomials which has nice properties. This point of view deserves to be examined in more detail. In particular, one interesting open problem would be to find closed (non-recursive) formulas for the coefficients $c_{n,m,k}$.

Note also that the recurrence \eqref{cnmkrecur} is not quite as arbitrary as it might appear at first
sight. The coefficients of $c_{n,m-1,k}$ on the right-hand side obey a symmetry relation under the interchange $k \leftrightarrow n+m+1-k$, namely
\bes
- 2k \leftrightarrow -2(n-k+m+1), \quad (n+4k-2m-2) \leftrightarrow (5n-4k+2m+2).
\ees
The significance of this observation is unclear to us at the moment.

\section{Boundary behavior: triple and quadruple refinements}
\label{sec:tripquadref}

We will now prove Theorems~\ref{thm:trip} and \ref{thm:quad}. The idea
is to equate the two formulas for the specialized partition function, one given by the definition in \eqref{defpartfn}, and the other by the determinant in \eqref{wrons}.

We now set the notation for use in the proofs. As mentioned before, we will consider
$f_{n}^{(m)}(u)$ to be a polynomial in $a(u)$ and $b(u)$. This polynomial turns out to be homogenous of
degree $3n-2$. Therefore, it makes sense to consider it as a polynomial in the single variable $t=a(u)/b(u)$. For convenience, we will use the
same notation $f_{n}^{(m)}(t)$ for the new function (where the derivative is still with respect to the variable $u$).
From Proposition~\ref{prop:fnmrecur}  and from \eqref{abidentities}, it follows easily that
\be \label{fnmt}
f_{n}^{(m)}(t) =
\left(\frac{\sqrt{3}}{2}\right)^{2n-1-m} 
b^{3n-2}(u) \left( t-1 \right)^{2n-1-m} \sum_{k=1}^{n+m} c_{n,m,k} t^{k-1}.
\ee
We will have at most four variables in the determinant we want to evaluate.
In continuation with Theorem~\ref{thm:wrons}, we use $u_{i}$ for the variables.
Throughout the computation, let $a_{i}$ (resp. $b_{i}$) be a shorthand for $a(u_{i})$ (resp. $b(u_{i})$), and denote $t_{i} = a_{i}/b_{i}$.

\subsection{Triple refinement}
\begin{proofof}{Theorem~\ref{thm:trip}}

We will use Theorem~\ref{thm:wrons} to calculate the number of ASMs where the
1 in the first row is in the $i$th column, the 1 in the first column is in the $j$th row, and the 1 in the last row is in the $k$th column. 
To relate the determinant in \eqref{wrons} to the generating function of the numbers $A_n^{\textrm{TLB}}(i,j,k)$, we need to calculate a determinant of a matrix of size 3 of derivatives $f_n^{(m)}(u_d)$, corresponding to $m=0,1,2$, expressing each of the derivatives in the form \eqref{fnmt}.
Let us compute these derivatives; we shall see that the result will involve the functions $\alpha_n(t), \beta_n(t), \gamma_n(t)$ from \eqref{defalph}. For $m=0$ we have
\begin{align*} 
f_{n}^{(0)}(t) &= \left(\frac{\sqrt{3}}{2}\right)^{2n-1} 
b^{3n-2}(u) \left( t-1 \right)^{2n-1} \sum_{k=1}^{n} A_{n,k} t^{k-1}
\\ &= \left(\frac{\sqrt{3}}{2}\right)^{2n-1} 
b^{3n-2}(u) \left( t-1 \right)^{2n-1} \alpha_{n}(t).
\end{align*}
Something interesting happens while computing $f_{n}^{(1)}(t)$. Applying the recurrence \eqref{cnmkrecur}, we get
\be
\begin{split}
f_{n}^{(1)}(t) =& 
\frac12 \left(\frac{\sqrt{3}}{2}\right)^{2n-2} 
b^{3n-2}(u)(t-1)^{2n-2} 
\\ & \times\sum_{k=1}^{n+1}\Big( 
- 2k A_{n,k+1}  + (n+4k-4) A_{n,k} 
\\ & \qquad \ \ \ \ \, + (5n-4k+4) A_{n,k-1} -2(n-k+2) A_{n,k-2} 
\Big) t^{k-1}.
\end{split}
\ee
Curiously, the coefficient of $t^{k-1}$ in this expression simplifies to a constant (depending on $n$ but not on $k$) times $A_{n-1,k-1}$, because of the identity
\begin{align*}
\frac{(n-2)!(3n-2)!}{(2n-3)!(2n-2)!} A_{n-1,k-1} =& 
- 2k A_{n,k+1}  + (n+4k-4) A_{n,k}  \\ &+ (5n-4k+4) A_{n,k-1} -2(n-k+2) A_{n,k-2}
\end{align*}
(which can be easily verified from \eqref{eq:refined-asm-thm} by elementary algebra). Thus we get that
\be \label{fnm1}
 f_n^{(1)}(t) = \left(\frac34\right)^{n-1}
\frac{(n-2)!(3n-2)!}{2(2n-3)!(2n-2)!}  b^{3n-2}(u)(t-1)^{2n-2} \beta_n(t).
\ee

Unfortunately, this nice phenomenon does not persist for higher derivatives. For the second derivative, using \eqref{cnmkrecur} again, after a short computation we get
\be \label{fnm2}
\begin{split}
f_{n}^{(2)}(t) =& \left(\frac{\sqrt{3}}{2}\right)^{2n-3}\frac{(n-2)!(3n-2)!}{4 (2n-3)!(2n-2)!}  b^{3n-2}(u) (t-1)^{2n-3} 
\\ & \times\sum_{k=1}^{n+2} \Big(
-2(n-k+3) A_{n-1,k-3} + (5n-4k+6) A_{n-1,k-2} \\& 
\qquad \ \ \ \  + (n+4k-6) A_{n-1,k-1} -2k A_{n-1,k}  \Big) \, t^{k-1} \\
=& 
\left(\frac{\sqrt{3}}{2}\right)^{2n-3}
\frac{(n-2)!(3n-2)!}{4(2n-3)!(2n-2)!}  b^{3n-2}(u)(t-1)^{2n-3} \gamma_{n}(t).
\end{split}
\ee
(The expression for $\gamma_n(t)$ does not seem to simplify further.)

Let us now look at the denominator on the right-hand side of  \eqref{wrons}. In the case $k=3$, it can be written (making use of \eqref{abidentities} again) as
\begin{align*}
\prod_{i=1}^{3} \sin^{2n-3}(u_{i}) &
\prod_{1 \leq i < j \leq 3} \!\! \sin(u_{i}-u_{j})
\\ &= \left(\frac34\right)^{3n-3} \prod_{i=1}^3 b_{i}^{2n-3}(t_{i}-1)^{2n-3}
\prod_{1 \leq i < j \leq 3} b_i b_j(t_i-t_j) 
\\ &= \left(\frac34\right)^{3n-3} \prod_{i=1}^3 b_i^{2n-1} (t_i-1)^{2n-3} \prod_{1\le i <j\le 3} (t_i-t_j).
\end{align*}
Combining this with the above computations, we find that in the case $k=3$, \eqref{wrons} transforms into
\be \label{zntriply1}
\begin{split}
Z_n&(u_1,u_2,u_3,0^{2n-3}) = \rho_n (b_1 b_2 b_3)^{n-1} \Delta(t_1,t_2,t_3)^{-1} 
\\ & \ \ \ \times \det\left( 
\begin{array}{ccc} (t_1-1)^2 \alpha_n(t_1) & (t_2-1)^2 \alpha_n(t_2) & (t_3-1)^2 \alpha_n(t_3) \\
(t_1-1) \beta_n(t_1) & (t_2-1) \beta_n(t_2) & (t_3-1) \beta_n(t_3) \\
\gamma_n(t_1) & \gamma_n(t_2) & \gamma_n(t_3)
\end{array}
\right),
\end{split}
\ee
where 
$\rho_n = \tfrac18 \zeta_{n,3} \left(\frac{(n-2)!(3n-2)!}{(2n-3)!(2n-2)!}\right)^2$ (which becomes \eqref{defrhon} upon substituting the value of $\zeta_{n,3}$ from \eqref{eq:formula-zetank}).

Now we move on to the left hand side of \eqref{wrons} with $k=3$. Our goal is to use the original definition of $Z_n(u_1,u_2,u_3,0^{2n-3})$ to represent it as a factor of $b_1 b_2 b_3$ times a generating function in $t_1, t_2, t_3$ expressed in terms of the coefficients $A_n^{\textrm{TLB}}$. Equating this expression to \eqref{zntriply1} will give the theorem.

First, note that since $Z_n$ is a symmetric function, we are free to choose which of the spectral parameters $x_1,\ldots,x_n,y_1,\ldots,y_n$ each of the variables $u_1, u_2, u_3$ is chosen to stand for. Since we are interested in relating $Z_n$ to the coefficients $A_n^{\textrm{TLB}}$, the natural choice is to set $x_1 = u_1$, $y_1=u_2$ and $x_n=u_3$; all other spectral parameters are of course set to $0$. For this choice, the weight $w(C)$ of each square ice configuration in the sum \eqref{defpartfn} is influenced only by the behavior of the configuration in the top and bottom rows and the leftmost column. In particular, we can divide the configurations into classes according to the position of the $1$s in the top and bottom rows and leftmost column in the corresponding alternating sign matrix. Within each class, all configurations have the same weight.

Let $C\in \mathcal{C}_n$ be a square ice configuration, and denote by $i,j,k$ the positions of the unique $1$ in the top row, leftmost column and bottom row, respectively, of the corresponding ASM.
There are three possible cases, each giving rise to a separate computation of the weight $w(C)$. 
\begin{enumerate}[1)]
\item $i=j=1$, $2\le k\le n$. The number of configurations of this type is $A_n^{\textrm{TLB}}(1,1,k)=A_{n-1,k-1}$, and the weight associated to such a configuration is given by
\be
\begin{split} \label{wt11}
w(C) &= b_1^{n-1} a_2^{n-2} a_3^{n-k} b_3^{k-2} b(u_3-u_2)
\\ &=
b_1^{n-1} a_2^{n-2} a_3^{n-k} b_3^{k-2} (a_2 a_3 + b_2 b_3 - b_2 a_3)
\\ &=
(b_1 b_2 b_3)^{n-1} (1-t_3+t_2 t_3) t_2^{n-2} t_3^{n-k}.
\end{split}
\ee
These weights arise as follows: in the first row, the first vertex has weight $c(u_1-u_2)=1$, and the remaining $n-1$ vertices have weight $b(u_1)=b_1$. In the first column, the vertices in positions $2$ through $n-1$ have weight $b(-u_2)=a(u_2)=a_2$, and the last vertex has weight $b(u_3-u_2)$, which can be rewritten as $a_2a_3+b_2b_3-b_2a_3$ (using \eqref{abidentities}). In the last row, the vertices in positions $2$ through $k-1$ have weight $b(u_3)=b_3$, the vertex in position $k$ has weight $c(u_3)=1$, and the vertices in positions $k+1$ through $n$ have weight $a(u_3)=a_3$.

\item $2\le i\le n$, $j=n$, $k=1$. The number of configurations is $A_{n-1,i-1}$, and the weight is seen (using similar reasoning to that explained above) to be
\be
\begin{split} \label{wtn1}
w(C) &= a_1^{i-2} b_1^{n-i} a_3^{n-1} b_2^{n-2} a(u_1-u_2)
\\ &=
a_1^{i-2} b_1^{n-i} a_3^{n-1} b_2^{n-2} (a_1 a_2 + b_1 b_2 -b_1 a_2)
\\ &=
(b_1 b_2 b_3)^{n-1} (1-t_2+t_1 t_2) t_1^{i-2} t_3^{n-1}.
\end{split}
\ee

\item $2\le i,k\le n$, $2\le j \le n-1$. In this case the number of configurations is $A_n^{\textrm{TLB}}(i,j,k)$. A short computation gives the weight of such a configuration as
\begin{align*}
w(C)&=
a_1^{i-2} b_1^{n-i} a_2^{n-j-1} b_2^{j-2} a_3^{n-k} b_3^{k-2}
a(u_1-u_2) b(u_3-u_2)
\\ &=
a_1^{i-2} b_1^{n-i} a_2^{n-j-1} b_2^{j-2} a_3^{n-k} b_3^{k-2}
\\ & \ \ \ \ \ \ \times(a_1 a_2 + b_1 b_2 - b_1 a_2) ( a_2 a_3 + b_2 b_3 - b_2 a_3)
\\ &=
(b_1 b_2 b_3)^{n-1} (1-t_2+t_1 t_2)(1-t_3+t_2 t_3) t_1^{i-2} t_2^{n-j-1} t_3^{n-k}.
\end{align*}
\end{enumerate}

Now all that is left to do is to combine the three cases to arrive at the partition function. We have
\begin{align*}
&\frac{1}{(b_1 b_2 b_3)^{n-1}}  Z_n(u_1,u_2,u_3,0^{2n-3})
=
(1-t_3+t_2 t_3) \sum_{k=2}^n A_{n-1,k-1} t_2^{n-2} t_3^{n-k}
\\ & \quad +
(1-t_2+t_1 t_2) \sum_{i=2}^n A_{n-1,i-1} t_1^{i-2} t_3^{n-1}
\\ & \quad +
(1-t_2+t_1 t_2)(1-t_3+t_2 t_3) \sum_{i=2}^n \sum_{j=2}^{n-1} \sum_{k=2}^n
A_n^{\textrm{TLB}}(i,j,k) t_1^{i-2} t_2^{n-j-1} t_3^{n-k}.
\end{align*}
Comparing this to the result \eqref{zntriply1} (and noting that the first two sums can be expressed in terms of $\alpha_{n-1}(\cdot)$, with a small simplification in the first sum that arises using the symmetry relation $A_{n-1,j}=A_{n-1,n-j}$) gives us just the equality of generating functions claimed in \eqref{tripref} (except with the variables $x,y,z$ replaced by $t_1,t_2,t_3$), and finishes the proof.
\end{proofof}

\subsection{Quadruple refinement}

\begin{proofof}{Theorem~\ref{thm:quad}}
The strategy is the same as that of Theorem~\ref{thm:trip}.
Recall that
\be
A_n^{\textrm{TLBR}}(i,j,k,\ell) = \#\Big\{ (m_{i,j})_{i,j=1}^n \in \asm_n\ |\ m_{1,i}\!=\!m_{j,1}\!=\!m_{n,k}\!=\!m_{\ell,n}\!=\!1 \Big\}.
\ee

As before, we will use Proposition~\ref{prop:fnmrecur} to calculate $f_{n}^{(3)}$.
Since the summand in the expansion of $f_{n}^{(2)}(t)$ in \eqref{fnm2} is quite complicated, we leave it to the reader to verify that the recursive formula indeed yield a function proportional to $\delta_{n}(t)$ \eqref{defalph}, 
\be \label{fnm3}
\begin{split}
f_{n}^{(3)}(t)&= 
\left(\frac{\sqrt{3}}{2}\right)^{2n-4} 
b^{3n-2}(u) \left( t-1 \right)^{2n-4} \sum_{k=1}^{n+3} c_{n,3,k} t^{k-1}, \\
&= \left(\frac 34\right)^{n-2} 
\frac{(n-2)!(3n-2)!}{8(2n-3)!(2n-2)!}
b^{3n-2}(u) (t-1)^{2n-4} \delta_{n}(t).
\end{split}
\ee

The denominator on the right-hand side of  \eqref{wrons} in the case $k=4$ can be written as
\begin{align*}
\prod_{i=1}^4 \sin^{2n-4}(u_{i}) &
\prod_{1 \leq i < j \leq 4} \!\! \sin(u_{i}-u_{j})
\\ &= \left(\frac34\right)^{4n-5} \prod_{i=1}^4 b_{i}^{2n-4}
(t_{i}-1)^{2n-4}
\prod_{1 \leq i < j \leq 4} b_i b_j(t_i-t_j) 
\\ &= \left(\frac34\right)^{4n-5} \prod_{i=1}^4 b_i^{2n-1} (t_i-1)^{2n-4} \prod_{1\le i <j\le 4} (t_i-t_j).
\end{align*}

Combining the calculation of $f_{n}^{(3)}(t)$ in \eqref{fnm3} 
with those in
\eqref{fnm1} and \eqref{fnm2}, we find that in the case $k=4$, \eqref{wrons} transforms into
\be \label{znquadly1}
\begin{split}
Z_n&(u_1,u_2,u_3,u_{4},0^{2n-4}) = \sigma_n 
(b_1 b_2 b_3 b_{4})^{n-1} \Delta(t_1,t_2,t_3,t_{4})^{-1} 
\\ & \ \ \ \times 
\det\left( \begin{array}{cccc} 
\scs  (t_1-1)^3 \alpha_n(t_1) & \scs(t_2-1)^3\alpha_n(t_2) & \scs(t_3-1)^3\alpha_n(t_3) & \scs(t_4-1)^3 \alpha_n(t_4) \\
\scs  (t_1-1)^2 \beta_n(t_1) & \scs(t_2-1)^2 \beta_n(t_2) & \scs(t_3-1)^2\beta_n(t_3) & \scs(t_4-1)^2 \beta_n(t_4) \\
\scs  (t_1-1) \gamma_n(t_1) & \scs(t_2-1) \gamma_n(t_2) & \scs(t_3-1)\gamma_n(t_3) & \scs(t_4-1) \gamma_n(t_4) \\
\scs  \delta_n(t_1) & \scs\delta_n(t_2) & \scs\delta_n(t_3) & \scs\delta_n(t_4)
\end{array}\right),
\end{split}
\ee
where 
$\sigma_n = \tfrac1{64} \zeta_{n,4} \left(\frac{(n-2)!(3n-2)!}{(2n-3)!(2n-2)!}\right)^3$, which is the same as the right-hand side of \eqref{defsigman}.

Now we compute the partition function 
$Z(u_{1},u_{2},u_{3},u_{4},0^{2n-4})$
from the definitions. This time we will set the spectral parameters as $x_{1}=u_{1}, y_{1}=u_{2}, x_{n}=u_{3}$ and 
$y_{n}= u_{4}$.
Let $C\in \mathcal{C}_n$ be a square ice configuration, and
 we will use the variables $i,j,k,\ell$ to denote the positions of the $c$-type vertices in the first row, leftmost column, last row and rightmost column, respectively. These correspond to the spectral parameters $u_{1}, u_{2}, u_{3}, 
u_{4}$. Remember that $i$ and $k$ are column parameters and $j$ and $l$ are row parameters. 
Now there are seven possible cases depending on the values of $i,j,k,l$.

\begin{enumerate}
\item $i=j=1, k=\ell=n$.
The number of such configurations is clearly $A_{n-2}$. The weight is given as follows: both the $(1,1)$ and the $(n,n)$ vertices get a weight of $c(u_{1}-u_{2})=1$ and $c(u_{3}-u_{4})=1$ respectively. All other boundary vertices are of $b$-type. 
The other two corner vertices at $(1,n)$ and $(n,1)$ get a weight 
of $b(u_{3}-u_{2})$ and $b(u_{1}-u_{4})$ respectively.
Therefore we have $n-2$ factors of each of $b(u_{1})=b_{1}, b(-u_{2})=a_{2}, b(u_{3})=b_{3}$ and $b(-u_{4})=a_{4}$. Thus the total contribution to $Z_n$ from such configurations is
\be
\begin{split}
&A_{n-2} (b_{1}a_{2}b_{3}a_{4})^{n-2} \; b(u_{3}-u_{2}) 
\; b(u_{1}-u_{4}), \\
&=A_{n-2} (b_{1}a_{2}b_{3}a_{4})^{n-2} 
(a_{2}a_{3}+b_{2}b_{3}- b_{2}a_{3})\; 
(a_{1}a_{4}+b_{1}b_{4}- b_{4}a_{1}), \\
&=A_{n-2} (b_{1}b_{2}b_{3}b_{4})^{n-1} 
\; (1-t_{3}+t_{2} t_{3}) (1-t_{1}+t_{4} t_{1}) \left(t_{2} t_{4} \right)^{n-2}.
\end{split}
\ee

\item $i=j=n, k=\ell=1$. The reasoning in this case is essentially identical to that of the case above. This time all boundary vertices except the ones at $(1,n)$ and $(n,1)$ are of 
$a$-type. We eventually obtain a total contribution of 
\be
A_{n-2} (b_{1}b_{2}b_{3}b_{4})^{n-1} 
\; (1-t_{2}+t_{1} t_{2}) (1-t_{4}+t_{3} t_{4}) \left(t_{1} t_{3} \right)^{n-2}.
\ee

\item $i=j=1, 1<k, \ell < n$.
The next four cases follow the same kind of reasoning as the one here. We describe this in detail and leave it to the reader to 
verify the next three cases by analogous reasoning.
The number of configurations of this type is 
$$A_n^{\textrm{TLBR}}(1,1,k,\ell)=A^{\textrm{BR}}_{n-1}(k-1,\ell-1),$$ 
which is the number of ASMs of size $n-1$ with a 1 in the bottom row at column
$k-1$ and a 1 in the rightmost column at row $\ell-1$.
The weight associated to such a configuration is given by
\be
\begin{split} 
w(C) =& b_1^{n-2} a_2^{n-2} a_3^{n-1-k} b_3^{k-2} a_{4}^{\ell-2}
b_{4}^{n-1-\ell} b(u_3-u_2) b(u_{1}-u_{4}) a(u_{3}-u_{4})
\\ =&
b_1^{n-2} a_2^{n-2} a_3^{n-1-k} b_3^{k-2} a_{4}^{\ell-2}
b_{4}^{n-1-\ell} (a_2 a_3 + b_2 b_3 - b_2 a_3)\\
& \times
(a_{1}a_{4}+b_{1}b_{4}- b_{4}a_{1})
(a_{3}a_{4}+b_{3}b_{4}- b_{3}a_{4})
\\ =&
(b_1 b_2 b_3 b_{4})^{n-1} t_2^{n-2} t_3^{n-1-k} 
t_{4}^{\ell-2}(1-t_3+t_2 t_3) \\
& \times
(1-t_{4}+t_{3} t_{4}) (1-t_{1}+t_{4} t_{1}).
\end{split}
\ee
These weights arise as follows: in the first row, the first vertex has weight $c(u_1-u_2)=1$, the next remaining $n-1$ vertices have weight $b(u_1)=b_1$ and the last vertex has weight 
$b(u_{1}-u_{4})$. In the first column, the vertices in positions $2$ through $n-1$ have weight $b(-u_2)=a(u_2)=a_2$, and the last vertex has weight $b(u_3-u_2)$, which can be rewritten as $a_2a_3+b_2b_3-b_2a_3$ (using \eqref{abidentities}). In the last row, the vertices in positions $2$ through $k-1$ have weight $b(u_3)=b_3$, the vertex in position $k$ has weight $c(u_3)=1$, the vertices in positions $k+1$ through $n-1$ have weight $a(u_3)=a_3$, and the last vertex has weight $a(u_{3}-u_{4})$. In the last column,  the vertices in positions $2$ through $\ell-1$ have weight $b(-u_4)=a_4$, the vertex in position $\ell$ has weight 
$c(u_4)=1$ and the vertices in positions $\ell+1$ through $n-1$ have weight $a(-u_4)=b_4$.

We now note that $A_{n-1}^{\textrm{BR}}(k-1,\ell-1) = 
A_{n-1}^{\textrm{TL}}(n+1-k,n+1-\ell)$ by rotational symmetry. 
Using the generating function defined in \eqref{deftopleftgf},
the sum over $k,\ell$ after elementary manipulations is given by 
\begin{align*}
&(b_1 b_2 b_3 b_{4})^{n-1} \; t_2^{n-2} \; t_{4}^{n-3}
(1-t_3+t_2 t_3) \\
&(1-t_{4}+t_{3} t_{4}) (1-t_{1}+t_{4} t_{1}) 
\mathcal A_{n-1}^{\textrm{TL}} \left(t_{3}, \frac 1{t_{4}} \right).
\end{align*}

\item $j=n,k=1, 1<i, \ell < n$.
 There are 
$A_{n-1}^{\textrm{TR}}(i-1,\ell) = 
A_{n-1}^{\textrm{TL}}(n+1-i,\ell)$ such configurations and the weight of such a configuration is given by 
\be
\begin{split}
w(C)=&(b_1 b_2 b_3 b_{4})^{n-1} 
t_{1}^{i-2} \; t_{4}^{\ell-2} t_{3}^{n-2} \\
& \times (1-t_{2}+t_{1} t_{2}) (1-t_{4}+t_{3} t_{4}) 
(1-t_{1}+t_{4} t_{1}) .
\end{split}
\ee
The sum over $i,\ell$ is now given by 
\begin{align*}
&(b_1 b_2 b_3 b_{4})^{n-1} \; t_3^{n-2} \; t_{1}^{n-3}
(1-t_2+t_1 t_2) \\
&(1-t_{4}+t_{3} t_{4}) (1-t_{1}+t_{4} t_{1}) 
\mathcal A_{n-1}^{\textrm{TL}} \left(t_{4}, \frac 1{t_{1}} \right).
\end{align*}

\item $k=\ell=n, 1<i, j < n$.
There are $A_{n-1}^{\textrm{TL}}(i,j)$ such configurations and the weight of such a configuration is given by 
\be
\begin{split}
w(C)=&(b_1 b_2 b_3 b_{4})^{n-1} t_{1}^{i-2} \; t_{2}^{n-1-j}
\; t_{4}^{n-2} \\
& \times (1-t_{2}+t_{1} t_{2}) (1-t_{3}+t_{2} t_{3})  
(1-t_{1}+t_{4} t_{1}) 
\end{split}
\ee
The sum over $i,j$ is given by 
\begin{align*}
&(b_1 b_2 b_3 b_{4})^{n-1} \; t_4^{n-2} \; t_{2}^{n-3}
(1-t_3+t_2 t_3) \\
&(1-t_{2}+t_{1} t_{2}) (1-t_{1}+t_{4} t_{1}) 
\mathcal A_{n-1}^{\textrm{TL}} \left(t_{1}, \frac 1{t_{2}} \right).
\end{align*}

\item $i=n, \ell=1, 1<j,k < n$.
There are $A_{n-1}^{\textrm{BL}}(k,j-1)=
A_{n-1}^{\textrm{TL}}(k,n+1-j)$ such configurations and the weight of such a configuration is given by 
\be
\begin{split}
w(C)=&(b_1 b_2 b_3 b_{4})^{n-1} t_{1}^{n-2}\; 
t_{2}^{n-1-j} \; t_{3}^{n-1-k}
 \\
&\times(1-t_{2}+t_{1} t_{2}) (1-t_{3}+t_{2} t_{3}) 
(1-t_{4}+t_{3} t_{4}) t_{1}^{n-2}.
\end{split}
\ee
The sum over $j,k$ is given by 
\begin{align*}
&(b_1 b_2 b_3 b_{4})^{n-1} \; t_1^{n-2} \; t_{3}^{n-3}
(1-t_3+t_2 t_3) \\
&(1-t_{4}+t_{3} t_{4}) (1-t_{2}+t_{1} t_{2}) 
\mathcal A_{n-1}^{\textrm{TL}} \left(t_{2}, \frac 1{t_{3}} \right).
\end{align*}

\item $1 < i,j,k,\ell < n$. 
The number of such matrices $A_n^{\textrm{TLBR}}(i,j,k,\ell)$ is the quantity we are ultimately interested in calculating. 
We describe in detail the weight of vertices in the first row.
The weights of vertices in the first column, last row and last 
column are obtained by similar arguments.
In the first row, the first vertex has weight $a(u_1-u_2)$, the next $i-2$ vertices have weight $a(u_1)=a_1$, the $i$th vertex has weight $c(u_{1})=1$, the next $n-i-1$ vertices have weight 
$b(u_{1})=b_{1}$, and the last vertex has weight $b(u_{1}-u_{4})$. The weight comes out to
\be
\begin{split}
w(C) =& \; a_{1}^{i-2} \; b_1^{n-i-1} \; a_2^{n-1-j} \; b_{2}^{j-2}
\; a_3^{n-1-k} \; b_3^{k-2} \; a_{4}^{\ell-2} \; b_{4}^{n-1-\ell} \\
& \times
a(u_1-u_2) b(u_3-u_2) b(u_{1}-u_{4}) a(u_{3}-u_{4})
\\ =&
\; a_{1}^{i-2} \; b_1^{n-i-1} \; a_2^{n-1-j} \; b_{2}^{j-2} 
\; a_3^{n-1-k} \; b_3^{k-2} \; a_{4}^{\ell-2} \; b_{4}^{n-1-\ell} 
\\
& \times
(a_2 a_3 + b_2 b_3 - b_2 a_3)
(a_{1}a_{4}+b_{1}b_{4}- b_{4}a_{1})\\
& \times
(a_{3}a_{4}+b_{3}b_{4}- b_{3}a_{4})
(a_{1}a_{2}+b_{1}b_{2}- b_{1}a_{2})
\\ =&
(b_1 b_2 b_3 b_{4})^{n-1} \; t_{1}^{i-2} \; t_2^{n-1-j} 
\; t_3^{n-1-k} t_{4}^{\ell-2}
 \\
& \times(1-t_{2}+t_{1} t_{2}) (1-t_{3}+t_{2} t_{3}) 
(1-t_{4}+t_{3} t_{4}) (1-t_{1}+t_{4} t_{1}).
\end{split}
\ee

\end{enumerate}

As in the triply-refined case, comparing the sum of the contributions from these seven cases
to the result \eqref{znquadly1} gives us the claimed identity~\eqref{quadref}.
\end{proofof}

\appendix

\section{Appendix: The Colomo-Pronko formula for the partition function and Theorem~\ref{thm:tripquad-simplified}}

Colomo and Pronko guessed in \cite{colomopronko-orthogonal} and proved
in \cite{colomopronko-emptiness} a determinantal formula for the
partition function $Z_n$ with $k$ inhomogeneous spectral parameters and $2n-k$ homogeneous ones (that is, set to $0$, in our notation), in the case $k\le n$.
After a small adaptation to translate the formula to our notation, their formula (Eq.~(5.8) in \cite{colomopronko-orthogonal}) reads
\be \label{eq:colomopronko-formula}
\begin{split}
Z_n&(u_1,\ldots,u_k,0^{2n-k}) \\ &= \Pi_n \frac{1}{\Delta(s_1,\ldots,s_k)}
\det\left( s_i^{k-j} (s_i-1)^{j-1} H_{n-k+j}(s_i) \right)_{i,j=1}^k,
\end{split}
\ee
where $H_n(u)$ is a single-variable generating function similar to our $\alpha_n(t)$ except that that it relates to a more general range of parameters of the square ice model, $s_i=s(u_i)$ is given by
\be
s_i = \frac{b(u_i)}{a(u_i)} = \frac{\sin\left(\frac{\eta}{2}-u_i\right)}{\sin\left(\frac{\eta}{2}+u_i\right)},
\ee
and $\Pi_n$ is a certain pre-factor which is easy to compute (see Eq.~(5.5) of \cite{colomopronko-orthogonal}).
This formula is very similar to \eqref{wrons} and some comments are in order:
\begin{enumerate}

\item The formula \eqref{eq:colomopronko-formula} works only for the case $k \leq n$, since Colomo and Pronko consider the case where all the row spectral variables are set to $0$. Our formula \eqref{wrons}
  works for all $k\le 2n$. On the other hand, the Colomo-Pronko formula holds for all
  values of the crossing parameter.

\item A comparison of \eqref{eq:colomopronko-formula} with \eqref{wrons} suggests that
  our function $f^{(m)}_n(t)$ can be expressed as a linear combination
  of $t^{j} (t-1)^{m-j} f_{n-j}(t)$ for $j$ from 0 to $m$. Numerical experiments
  suggest that this is always true.
  However, we are not aware of
  an easy way to see this.
  
\end{enumerate}

Motivated by Colomo's observation, we used computer algebra 
to verify that the modified function $\tilde{\gamma}_n(t)$ defined in \eqref{redefalph} can be represented as a linear combination of $\alpha_n(t), \beta_n(t), \gamma_n(t)$ from \eqref{defalph}, and similarly the modified $\tilde{\delta}_n(t)$ can be represented as a linear combination of the functions $\alpha_n(t), \beta_n(t), \gamma_n(t), \delta_n(t)$. If we denote by $\mu_n$ and $\nu_n$, the coefficient of $\gamma_n(t)$ in the first linear combination and the coefficient of $\delta_n(t)$ in the second one, respectively, then this fact
establishes rigorously the claim of Theorem~\ref{thm:tripquad-simplified}.
The computation is straightforward with the help of a symbolic algebra package, and is included in a \texttt{Mathematica} notebook \texttt{GammaDelta} accompanying this paper, which the interested reader may download from the authors' websites.

\bigskip
\noindent
\textbf{Authors' contact information} \\
Department of Mathematics \\ University of California, Davis \\ One Shields Ave. \\ Davis, CA 95616 \\ USA
\\ \ \\
Email: \texttt{ayyer@math.ucdavis.edu}, \texttt{romik@math.ucdavis.edu}

\end{document}